\documentclass[11pt]{amsart}
\usepackage{mathrsfs}
\usepackage{}
\usepackage{amsfonts}
\usepackage{amssymb}
\usepackage{bbm}

\newtheorem{thm}{Theorem}[section]
\newtheorem{cor}[thm]{Corollary}
\newtheorem{lem}[thm]{Lemma}

\newtheorem{prop}[thm]{Proposition}
\newtheorem{rem}[thm]{Remark}
\numberwithin{equation}{section}

\newcommand{\be}{\begin{equation}}
\newcommand{\ee}{\end{equation}}
\newcommand{\bes}{\begin{eqnarray}}
\newcommand{\ees}{\end{eqnarray}}
\newcommand{\bess}{\begin{eqnarray*}}
\newcommand{\eess}{\end{eqnarray*}}

\newcommand{\bali}{\begin{align}}
\newcommand{\eali}{\end{align}}

\def\kk{\mathbbm{k}}
\def\Hd{H_{\mathcal {D}}}

\begin{document}
\title[Green Rings of Pointed Rank One Hopf algebras]{Green Rings of Pointed Rank One Hopf algebras of Non-nilpotent Type}
\author{Zhihua Wang}
\address{Z. Wang\newline School of Mathematical Science, Yangzhou University,
Yangzhou 225002, China}
\email{mailzhihua@gmail.com}
\author{Libin Li}
\address{L. Li\newline School of Mathematical Science, Yangzhou University,
Yangzhou 225002, China}
\email{lbli@yzu.edu.cn}
\author{Yinhuo Zhang}
\address{Y. Zhang\newline Department of Mathematics and Statistics, University of Hasselt, Universitaire Campus, 3590 Diepeenbeek, Belgium }
\email{yinhuo.zhang@uhasselt.be}
\date{}

\begin{abstract} In this paper, we continue our study of the Green rings of finite dimensional pointed Hopf algebras of rank one initiated in \cite{WLZ}, but focus on those Hopf algebras of non-nilpotent type. Let $H$ be a finite dimensional pointed rank one Hopf algebra of non-nilpotent type. We first determine all non-isomorphic indecomposable $H$-modules and describe the Clebsch-Gordan formulas for them. We then study the structures of both the Green ring $r(H)$ and the Grothendieck ring $G_0(H)$ of $H$ and establish the precise relation between the two rings. We use the Cartan map of $H$ to study the Jacobson radical and the idempotents of $r(H)$. It turns out that the Jacobson radical of $r(H)$ is exactly the kernel of the Cartan map, a principal ideal of $r(H)$, and $r(H)$ has no non-trivial idempotents. Besides, we show that the stable Green ring of $H$ is a transitive fusion ring. This enables us to calculate Frobenius-Perron dimensions of objects of the stable category of $H$. Finally, as an example, we present both the Green ring and the Grothendieck ring of the Radford Hopf algebra in terms of generators and relations.
\end{abstract}
\maketitle

\section{Introduction }
In \cite{WLZ}, we studied the Green rings of finite dimensional pointed rank one Hopf algebras of nilpotent type. The Green ring of a finite dimensional pointed rank one Hopf algebra $H$ of nilpotent type is isomorphic to a polynomial ring in one variable over the Grothendieck ring of $H$ modulo one relation. This Green ring is a symmetric algebra over $\mathbb{Z}$ with the Jacobson radical being a principal ideal. The stable Green ring (i.e., the Green ring of the stable category) of $H$ is isomorphic to the quotient ring of the Green ring of $H$ modulo the ideal generated by isomorphism classes of indecomposable projective $H$-modules. In particular, the complexified stable Green algebra is both a group-like algebra and a bi-Frobenius algebra introduced by Doi and Takeuchi (cf. \cite{Doi1, Doi2, Doi3, DT}).

In this paper, we continue the study of the Green rings of finite dimensional pointed Hopf algebras of rank one, but concentrate us on pointed rank one Hopf algebras of non-nilpotent type. The representation theory of a pointed rank one Hopf algebra of non-nilpotent type is closely related to that of some pointed rank one Hopf algebra of nilpotent type. However, as we shall see in the non-nilpotent case, the Green ring is much more complicated.  The main reason is due to the lack of the Chevalley property
in the non-nilpotent case.  This leads to the problem that the free abelian group generated by simple modules doest not form a subring of the Green ring. Because of this, it is hard to describe the Green ring in terms of generators and relations as we did in the nilpotent case \cite{WLZ}. However, we can still study the ring-theoretical properties of the Green ring, including the Frobenius property, the Jacobson radical and idempotents of the Green ring.

The paper is organized as follows. In Section 2 we first recall the Hopf algebra structure of a finite dimensional pointed rank one Hopf algebra $H$ of non-nilpotent type. We construct a central idempotent element $e$ of $H$ such that the subalgebra $H(1-e)$ is  a Hopf ideal of $H$. It turns out  that the quotient Hopf algebra $\overline{H}:=H/H(1-e)$ is a pointed rank one Hopf algebra of nilpotent type. Thus, the indecomposable $\overline{H}$-modules are also indecomposable modules over $H$.
On the other hand, every finite dimensional simple module over the subalgebra $H(1-e)$ is projective. Consequently,  $H(1-e)$ is a semisimple subalgebra of $H$. As a result, we obtain all indecomposable $H$-modules from both the indecomposable modules over $\overline{H}$ and the simple modules over $H(1-e)$.

In Section 3 and Section 4, we study the Clebsch-Gordan formulas
for the decompositions of tensor products of  indecomposable $H$-modules. This includes the decompositions of the tensor products of simple $H(1-e)$-modules. We work out the relation between the Green rings $r(H)$ of $H$ and $r(\overline{H})$ of $\overline{H}$. Let $\mathcal{P}$ be the ideal of $r(H)$ generated by isomorphism classes of indecomposable projective $H$-modules and the direct sum $r(\overline{H})\oplus \mathcal{P}$ a ring with a new multiplication rule. Instead of describing the Green ring $r(H)$ in terms of generators and relations, we show that $r(H)$ is isomorphic to a quotient ring of $r(\overline{H})\oplus \mathcal{P}$ modulo an ideal described in detail.

In Section 5, we investigate the Grothendieck ring $G_0(H)$ of $H$. We first show that the Green ring $r(H)$ admits an associative, symmetric and nonsingular $\mathbb{Z}$-bilinear form. The Grothendieck ring $G_0(H)$ is isomorphic to the quotient ring $r(H)/\mathcal {P}^{\perp}$, where
$\mathcal {P}^{\perp}$ is the orthogonal complement of $\mathcal {P}$ with respect to the bilinear form.
In addition, $G_0(H)$ can also be regarded as a subring of the Green ring $r(\kk G)$ of the group algebra $\kk G$, where $G$ is the group of group-like elements of $H$. Besides, we describe the Cartan matrix of the Cartan map from $\mathcal {P}$ to $G_0(H)$ with respect to the given bases. It turns out that this Cartan matrix is of block-diagonal form with each block described explicitly.

In section 6, we use the Cartan map to study the Jacobson radical and idempotents of $r(H)$.
Our approach to the Jacobson radical of $r(H)$ is different from the nilpotent case carried out in \cite[Section 5]{WLZ}, where we had to study the irreducible representations of the complexified Green ring $r(H)\otimes_{\mathbb{Z}}\mathbb{C}$ in order to get the rank of the Jacobson radical of the Green ring $r(H)$. With the new approach of using the Cartan map,  we can show that the Jacobson radical of $r(H)$ is precisely the kernel of the Cartan map, which is exactly the intersection $\mathcal {P}\cap\mathcal {P}^{\perp}$. This implies that the Jacobson radical of $r(H)$ is a principal ideal and that the Green ring $r(H)$ has no non-trivial idempotents, a property that holds as well in the nilpotent case.

In section 7, we consider the stable Green ring of $H$. Based on the representations of $H$, we see that the stable Green rings $r_{st}(H)$ of $H$ and $r_{st}(\overline{H})$ of $\overline{H}$ coincide, where the latter has been studied in \cite[Section 6]{WLZ}. Here we go further to show that the stable Green ring $r_{st}(H)$ is a transitive fusion ring. This allow us to study the Frobenius-Perron dimensions of objects in the stable category of $H$ in the framework of $r_{st}(H)$. Using the Dickson polynomials (of the second type), we calculate all Frobenius-Perron dimensions of indecomposable objects of the stable category of $H$.

In Section 8, as an example, we compute the Green ring $r(H)$ of the Radford Hopf algebra $H$. Although the Green ring of a pointed rank one Hopf algebra  of non-nilpotent type is difficult to describe in terms of generators and relations, we can do so for the Radford Hopf algebra $H$ bacause the quotient Hopf algebra $\overline{H}=H/H(1-e)$ is a Taft Hopf algebra. According to \cite[Theorem 3.10]{Chen} or \cite[Theorem 4.3]{WLZ}, we are able to show that the Green ring $r(H)$ is isomorphic to a quotient ring of a polynomial ring in several variables over $\mathbb{Z}$ modulo certain relations described explicitly. As a quotient of the Green ring $r(H)$, the Grothendieck ring $G_0(H)$ of $H$ can  be presented as well in terms of generators and relations.

Throughout, we work over an algebraically closed field $\mathbbm{k}$
of characteristic zero. Unless otherwise stated, all algebras, Hopf
algebras and modules are defined over $\mathbbm{k}$; all modules are
left modules and finite dimensional; $\dim$ and $\otimes$ stand for
$\dim_{\mathbbm{k}}$ and $\otimes_{\mathbbm{k}}$ respectively. For a finite set $B$, $|B|$ is the cardinality of $B$ and sp$B$ is the $\kk$-vector space spanned by $B$. We refer to \cite{Kas, Mon} for the theory of Hopf algebras.

\section{Indecomposable Representations}
\setcounter{equation}{0}

A quadruple $\mathcal {D}=(G,\chi,g,\mu)$ is
called a \emph{group datum} if $G$ is a finite group, $\chi$  a $\mathbbm{k}$-linear character of $G$, $g$ an element in
the center of $G$,
and $\mu\in\mathbbm{k}$ subject to $\chi^n=1$ or $\mu(g^n-1)=0$,
where $n$ is the order of $\chi(g)$. The
group datum $\mathcal {D}$ is said to be of \emph{nilpotent type} if $\mu(g^n-1)=0$, and it is of \emph{non-nilpotent type} if $\mu(g^n-1)\neq0$ and $\chi^n=1$ (cf. \cite{CHYZ, KR}).

Let $\Hd$ be the Hopf algebra associated to the group datum $\mathcal {D}=(G,\chi,g,\mu)$. Namely, as an associative algebra, $\Hd$ is generated by $y$ and all elements $h$ in $G$ such that $\kk G$ is a subalgebra of $\Hd$ and
\begin{equation}\label{equT10}y^n=\mu(g^n-1),\ yh=\chi(h)hy,\ \text{for}\ h\in G.\end{equation}
The comultiplication $\bigtriangleup$, the counit
$\varepsilon$, and the antipode $S$ of $\Hd$ are given respectively by
\begin{equation}\label{3.2}\bigtriangleup(y)=y\otimes g+1\otimes y,\ \varepsilon(y)=0,\ S(y)=-yg^{-1},\end{equation}
\begin{equation}\label{3.1}\bigtriangleup(h)=h\otimes h,\ \varepsilon(h)=1,\ S(h)=h^{-1},\end{equation}
for all $h\in G.$

The Hopf algebra $\Hd$ is a finite dimensional pointed Hopf algebra of rank one \cite{KR} with a canonical
$\mathbbm{k}$-basis $\{y^ih\mid h\in G,\ 0\leq i\leq n-1\}$. The group of group-like elements of
$H_{\mathcal {D}}$ is $G$ and $\dim H_{\mathcal {D}}=n|G|$. The Hopf algebra $\Hd$ is said to be of \emph{nilpotent type} (resp. \emph{non-nilpotent type}) if the associated group datum $\mathcal {D}$ is nilpotent (resp. non-nilpotent), see \cite{KR, WLZ}.

Throughout this paper, $\mathcal {D}=(G,\chi,g,\mu)$ is a fixed group datum of non-nilpotent type with $n$ the order of $q:=\chi(g)$. Without loss of generality we may assume that $\mu=1$ (see \cite[Corollary 1]{KR}). For the simplicity, we shall drop $\mathcal{D}$ from $\Hd$, and denote by $H$ the Hopf algebra $\Hd$ associated to the group datum $\mathcal{D}=(G,\chi, g, 1)$. In this case, the relations in (\ref{equT10}) become $$y^n=g^n-1,\ yh=\chi(h)hy,\ \textrm{for}\ h\in G.$$
Note that $\chi^n=1$ and the order of $q=\chi(g)$ is $n$. The order of $\chi$ is also $n$. If $n=1$, the Hopf algebra $H$ is nothing but $\kk G$. To avoid this, we always assume that $n\geq2$. In this case, $\chi(g)\neq 1$,  implying that $g\neq 1$ and $\chi\neq\varepsilon$.

Note that the order of $q=\chi(g)$ is $n$. The order of $g$ in the group $G$ is $nr$, where
$r>1$ since the group datum $\mathcal {D}$ is of non-nilpotent type. Denote by $$N=\{1,g^n,g^{2n},\cdots,g^{(r-1)n}\}$$ and $$e=\frac{1}{r}(1+g^n+g^{2n}+\cdots+g^{(r-1)n}).$$
Then $N$ is a normal subgroup of $G$ and $e$ is a central
idempotent of $H$ satisfying $eg^n=e$. Let $\overline{G}=G/N$. The following lemma is obvious.

\begin{lem}\label{lemT2}
Let $\kk G(g^n-1)$ be the ideal of $\mathbbm{k}G$ generated by $g^n-1$.
Then $\mathbbm{k}G/\kk G(g^n-1)\cong\mathbbm{k}\overline{G}$.
\end{lem}

\begin{rem} Note that $1-e$ and $1-g^n$ are both central elements of $H$ and $1-e=(1-g^n)\frac{1}{r}\sum_{k=1}^{r-1}(1+g^n+\cdots+g^{(k-1)n})$ whereas $1-g^n=(1-e)(1-g^n).$
\begin{enumerate}
  \item[(1)] The ideals of $\kk G$ generated respectively by $1-e$ and $1-g^n$ coincide. It follows from Lemma \ref{lemT2} that $\kk \overline{G}\cong\kk Ge$.
  \item[(2)] The ideals of $H$ generated respectively by $1-e$ and $1-g^n$ coincide. Thus, $H(1-e)$ is a Hopf ideal of $H$ because $H(1-g^n)$ is.
\end{enumerate}
\end{rem}

Let $\pi$ be the natural projection from $G$
to $\overline{G}$ and $\overline{h}:=\pi(h)$ for $h\in G$. The character $\chi$ of $G$ induces a character $\overline{\chi}$ of $\overline{G}$
such that $\overline{\chi}\circ\pi=\chi.$  Denote by $\overline{\mathcal{D}}=(\overline{G},\overline{\chi},\overline{g},0)$. Then the group datum $\overline{\mathcal{D}}$ is of nilpotent type since $\overline{g}^n-1=0$, where $n$ is the order of $\chi(g)=\overline{\chi}(\overline{g})$. Let $\overline{H}$ be the Hopf algebra associated to the group datum $\overline{\mathcal {D}}$.
More precisely,
$\overline{H}$ is generated as an algebra by $z$ and all
elements $\overline{h}\in \overline{G}$ such that $\kk\overline{G}$ is a subalgebra of $\overline{H}$ and the relations (\ref{equT10}), (\ref{3.2}) and (\ref{3.1}) are satisfied when we replace $y$ and elements $h\in G$ with $z$ and $\overline{h}\in\overline{G}$ respectively.
So $\overline{H}$ is a finite dimensional pointed rank one Hopf algebra of nilpotent type with a $\mathbbm{k}$-basis $\{z^i\overline{h}\mid\overline{h}\in \overline{G},0\leq i\leq n-1\}$.

The relation between Hopf algebras $H$ of non-nilpotent type and $\overline{H}$ of nilpotent type is described as follows:

\begin{prop}\label{propT2}
\begin{enumerate}
\item[$(1)$] $\overline{H}$ is isomorphic  to $H/H(1-e)$ as a Hopf algebra.
\item[$(2)$] $\overline{H}$ is isomorphic  to $He$ as an algebra.
\end{enumerate}
\end{prop}
\proof (1) It is easy to check that the algebra epimorphism
$\rho:H\rightarrow \overline{H}$ given by $\rho(y)=z$ and $\rho(h)=\overline{h}$, for any $h\in G$ respects the Hopf algebra
structure. The inclusion $H(1-e)\subseteq\ker\rho$ is obvious since $\rho(e)=1$. To verify that $\ker\rho=H(1-e)$, we only need to verify that the restriction of $\rho$ to the summand $He$ of $H$ is injective. In fact, if $\rho(\sum_{i=0}^{n-1}y^ia_ie)=0$ for $a_i\in \mathbbm{k}G$, then $\sum_{i=0}^{n-1}z^i\overline{a_ie}=0$, and hence each $\overline{a_ie}=0$. By Lemma \ref{lemT2}, $\overline{a_ie}=0$ if and only if $a_ie\in\kk G(g^n-1)$ if and only if $a_ie\in\kk G(1-e)$. It follows that $a_ie=0$, as desired.

(2) Follows from the algebra decomposition  $H=He\oplus H(1-e)$.\qed

Let $V$ be a $\kk G$-module and $x$ a formal variable. For any $k\in \mathbb{N}$, the set $x^kV=\{x^kv\mid v\in V\}$ is a $\kk$-vector space defined by $x^ku+x^kv=x^k(u+v)$ and $\lambda(x^kv)=x^k(\lambda v)$ for $u,v\in V$ and $\lambda\in\kk$. Obviously, $\dim x^kV=\dim V$. Moreover, $x^kV$ is  a $\kk G$-module with the $G$-action given by
\begin{equation}\label{nonnil}h(x^kv)=\chi^{-k}(h)x^khv,\end{equation} for any $h\in G$ and $v\in V$.

Let $\{V_i\mid i\in\Omega\}$ be a complete set of simple $\mathbbm{k}G$-modules up to isomorphism, where $V_0$ is the trivial module $\kk$. Since $g^n$ is a central element of $G$, the action of $g^n$ on each $V_i$ is a scalar multiple by a non-zero element of $\kk$, say $\lambda_i$. Let $\Omega_0=\{i\in\Omega\mid\lambda_i=1\}$ and $\Omega_1=\{i\in\Omega\mid\lambda_i\neq1\}$. In particular, $0\in\Omega_0$. It follows from Lemma \ref{lemT2} that $\{V_i\mid i\in\Omega_0\}$ is a complete set of simple $\mathbbm{k}\overline{G}$-modules up to isomorphism.

Let
$M(j,i):=V_i\oplus xV_i\oplus\cdots\oplus x^{j-1}V_i$ for $i\in\Omega_0$ and $1\leq j\leq n$.
Then $M(j,i)$ is an $H$-module, where the action of $h\in G$ on $M(j,i)$ is given by (\ref{nonnil}) and the action of $y$ on $M(j,i)$ is
$$y(x^kv)=\begin{cases}
x^{k+1}v, & 0\leq k\leq j-2,\\
0, & k=j-1,
\end{cases}$$
for any $v\in V_i$. Note that $\lambda_i=1$ for $i\in\Omega_0$. We have the following action of the idempotent element $e$:
\begin{equation}\label{equT6}
e(x^kv)=\frac{1}{r}\sum_{s=0}^{r-1}g^{sn}(x^kv)=x^k(\frac{1}{r}\sum_{s=0}^{r-1}g^{sn}v)=x^kv,
\end{equation}
for any $x^kv\in M(j,i)$. It follows from Proposition \ref{propT2} that each $M(j,i)$ is an $\overline{H}$-module. In particular, by \cite[Theorem 2.5]{WLZ}, $$\{M(j,i)\mid i\in\Omega_0,1\leq j\leq n\}$$ forms a complete set of finite dimensional indecomposable $\overline{H}$-modules up to isomorphism, where $M(1,0)\cong\kk$, $M(1,i)\cong V_i$ and $M(n,i)$ is the projective cover of $V_i$.

So far we have obtained all finite dimensional indecomposable modules over $\overline{H}\cong He$. Now we  determine all finite dimensional indecomposable modules over $H(1-e)$. It turns out that the indecomposable modules over $H(1-e)$ are induced from the simple $\kk G$-modules $V_j$, for $j\in \Omega_1$.

Let $P_j:=V_j\oplus xV_j\oplus\cdots\oplus x^{n-1}V_j$ for $j\in\Omega_1$.
Then $P_j$ is an $H$-module, where the action of $h\in G$ on $P_j$ is given by (\ref{nonnil}) and the element $y$ acts on $P_j$ as follows:
$$y(x^kv)=\begin{cases}
x^{k+1}v, & 0\leq k\leq n-2,\\
(\lambda_j-1)v, & k=n-1,
\end{cases}$$
for any $v\in V_j$. Note that $\lambda_j\neq1$ and $\lambda_j^r=1$ for $j\in\Omega_1$. Now the action of $e$ on $P_j$ is zero:
\begin{equation} \label{equT5}
e(x^kv)=\frac{1}{r}\sum_{s=0}^{r-1}g^{sn}(x^kv)=x^k(\frac{1}{r}\sum_{s=0}^{r-1}g^{sn}v)=x^k(\frac{1}{r}\sum_{s=0}^{r-1}\lambda_j^{s}v)=0,
\end{equation}
for any $x^kv\in P_j$. This means that each $P_j$ is in fact an $H(1-e)$-module.

For any $H$-module $V$, the subspace $V_y=\{v\in V\mid yv=0\}$ is a submodule of $V$. If $V_y=V$, then $V$ is called \emph{$y$-torsion}. If $V_y=0$, then $V$
is called \emph{$y$-torsionfree}. Obviously, if $V$ is a simple $H$-module, then $V$ is either $y$-torsion or $y$-torsionfree. It follows from Lemma \ref{lemT2} that an $H$-module $V$ is simple $y$-torsion if
and only if $V$ is a simple $\kk\overline{G}$-module. We have the following description of the simple $y$-torsionfree $H$-modules.

\begin{prop}\label{propT5}
\begin{enumerate}
\item[(1)] Every simple $y$-torsionfree $H$-module is isomorphic to $P_j$ for some $j\in\Omega_1$.
\item[(2)] Every $H$-module $P_j$, $j\in\Omega_1$,  is simple, projective and $y$-torsionfree.
\end{enumerate}
\end{prop}
\proof (1) If $V$ is a simple $y$-torsionfree $H$-module, then $V$ is semisimple as a $\mathbbm{k}G$-module
since $\mathbbm{k}G$ is a semisimple subalgebra of $H$. Hence there is some $j\in\Omega$ such that $V_j$ is a direct summand of $V$. By the simplicity of $V$, we have
$V=V_j+yV_j+\cdots+y^{n-1}V_j$ .
We claim that $j\in\Omega_1$ (i.e., $\lambda_j\neq1$) and
$V_j+yV_j+\cdots+y^{n-1}V_j$ is a
direct sum. In fact, if $\lambda_j=1$, then for any $v\in V_j$,
$y^nv=(g^n-1)v=(\lambda_j-1)v=0$. It follows that $v=0$ since $V$ is
$y$-torsionfree, a contradiction. Note that $g$ is a central element of $G$. There is a non-zero scalar $\omega_j$ such that $gv=\omega_jv$, for any $v\in V_j$. If $v_0+yv_1+\cdots+y^{n-1}v_{n-1}=0$ for $v_i\in V_j$,
then
\begin{align*}
0&=g^i(v_0+yv_1+\cdots+y^{n-1}v_{n-1})\\
&=\omega_j^i(v_0+q^{-i}yv_1+\cdots+q^{-(n-1)i}y^{n-1}v_{n-1}).
\end{align*}
Thus,
\begin{equation}\label{equT11}v_0+q^{-i}yv_1+\cdots+q^{-(n-1)i}y^{n-1}v_{n-1}=0,\ \textrm{for}\ 0\leq i\leq n-1.\end{equation}
This implies that $v_0=yv_1=\cdots=y^{n-1}v_{n-1}=0$ since the
coefficient matrix determined by Equations (\ref{equT11}) is a Vandermonde matrix (which is invertible). Now we have $V=V_j\oplus yV_j\oplus\cdots\oplus y^{n-1}V_j$ for some $j\in\Omega_1$, and hence $V$ is isomorphic to
$P_j$ with the isomorphism given by $$V\rightarrow P_j,\ \sum_{i=0}^{n-1}y^iv_i\mapsto\sum_{i=0}^{n-1}x^iv_i,$$  as desired.

(2) We first claim that $P_j=V_j\oplus xV_j\oplus\cdots\oplus
x^{n-1}V_j$ is $y$-torsionfree, for any $j\in\Omega_1$. Indeed, if  $y(v_0+xv_1+\cdots+x^{n-1}v_{n-1})=0$, then $$(g^n-1)(v_0+xv_1+\cdots+x^{n-1}v_{n-1})=y^n(v_0+xv_1+\cdots+x^{n-1}v_{n-1})=0.$$ It follows that $(\lambda_j-1)(v_0+xv_1+\cdots+x^{n-1}v_{n-1})=0$. We obtain that $v_0+xv_1+\cdots+x^{n-1}v_{n-1}=0$ since $\lambda_j\neq1$.

Now we check that $P_j$ is simple for any $j\in\Omega_1$. If $V$ is a non-zero simple submodule of $P_j$, then $V$ is simple $y$-torsionfree. By Part (1), $V\cong P_{j'}$ for some $j'\in\Omega_1$. Moreover, $V_{j'}\subseteq P_{j'}\cong V\subseteq P_j$. Note that $V_{j'}$ is a simple $\mathbbm{k}G$-module and $P_j$ is a direct sum with the summands of simple $\mathbbm{k}G$-modules. There exists some $k$ such that $V_{j'}\cong x^kV_j$. Thus, $\dim V_{j'}=\dim x^kV_j=\dim V_j$. This implies that $\dim V=\dim P_{j'}=\dim P_j$. As a result, $V=P_j$.

To see that $P_j$ is projective for any $j\in\Omega_1$, we denote by $e_j$ the primitive idempotent
of $\mathbbm{k}G$ such that $V_j\cong\mathbbm{k}Ge_j$ as $\mathbbm{k}G$-modules. Then $e_j$ is an idempotent of $H$ as well. Let $He_j$ be the left ideal of $H$ generated by $e_j$.  We have the following decomposition of $He_j$ into direct sum of simple $\mathbbm{k}G$-modules:
 $$He_j=\mathbbm{k}Ge_j\oplus
y\mathbbm{k}Ge_j\oplus\cdots\oplus y^{n-1}\mathbbm{k}Ge_j.$$
 Denote by $\zeta_j$ the isomorphism from $V_j$ to $\mathbbm{k}Ge_j$ and consider the
$\mathbbm{k}$-linear map as follows:
$$P_j\rightarrow He_j,\ \sum_{k=0}^{n-1}x^kv_k\mapsto\sum_{k=0}^{n-1}y^k\zeta_j(v_k).$$
It is straightforward to verify that this map is an $H$-module isomorphism.
\qed

We are now ready to determine all finite dimensional indecomposable modules over the subalgebra $H(1-e)$ of $H$.

\begin{thm}\label{3.4} The following hold:
\begin{enumerate}
\item[(1)] For any $j\in\Omega_1$, $P_j$ is a simple projective $H(1-e)$-module.
\item[(2)] Every simple $H(1-e)$-module is isomorphic to some $P_j$, $j\in\Omega_1$.
\item[(3)] The subalgebra $H(1-e)$ of $H$ is semisimple.
\item[(4)] The Jacobson radical of $H$ is a principal ideal generated by the element $ye$.
\end{enumerate}
\end{thm}
\proof (1) Follows from Proposition \ref{propT5} (2) and Equation (\ref{equT5}).

(2) If $V$ is a simple $H(1-e)$-module, then $V$ is a simple $H$-module via the projection from $H$ to  $H(1-e)$. Moreover, $V$ is  $y$-torsionfree. Indeed, if $yv=0$ for some $v\in V$, then
$(g^n-1)v=y^nv=0$ and hence $g^nv=v$. This implies that
$0=ev=\frac{1}{r}\sum_{k=0}^{r-1}g^{kn}v=v$. Thus, $V$ is a simple $y$-torsionfree $H$-module. It follows from
Proposition \ref{propT5} (1) that $V\cong P_j$ for some $j\in\Omega_1$.

(3) Follows from Part (1) and Part (2).

(4) Follows from the fact that the ideal $(ye)$ of $H$ generated by $ye$ is nilpotent and the quotient algebra $H/(ye)\cong H(1-e)\bigoplus\kk Ge$ is semisimple. \qed

\begin{rem}
Recall that a finite dimensional Hopf algebra over $\kk$ is said to have the \emph{Chevalley property} if the tensor product of any two simple modules is semisimple \cite[Definition 7.2.1]{Ge1}. One of the equivalent conditions is that the Jacobson radical of the Hopf algebra is a Hopf ideal \cite[Proposition 7.2.2]{Ge1}. Accordingly, the finite dimensional pointed rank one Hopf algebra $H$ of non-nilpotent type does not have the Chevalley property since the Jacobson radical of $H$ is not a Hopf ideal. As we shall see that the tensor product of any two simple $H$-modules is not necessary semisimple.
\end{rem}

In order to determine all simple $H(1-e)$-modules up to isomorphism, we need a permutation on the index set $\Omega$ of simple $\kk G$-modules.
Let $V_{\chi}$ and $V_{\chi^{-1}}$ be the two 1-dimensional simple
$\mathbbm{k}G$-modules corresponding to the $\kk$-linear character $\chi$ and $\chi^{-1}$ respectively. Similar to the nilpotent case in \cite{WLZ}, we have a unique permutation $\tau$ of the index set $\Omega$ determined by
$$V_{\chi^{-1}}\otimes V_s\cong V_s\otimes V_{\chi^{-1}}\cong V_{\tau(s)}$$ as $\kk G$-modules, for some $\tau(s)\in\Omega.$
Moreover, it is easy to see that
$s\in\Omega_0$ (resp. $s\in\Omega_1$) if and only if $\tau(s)\in\Omega_0$
(resp. $\tau(s)\in\Omega_1$), i.e., $\tau$ permutes the index set
$\Omega_0$ and $\Omega_1$ respectively.

\begin{lem}\label{lemT3}
For any $s\in\Omega$ and $t\in\mathbb{Z}$, the following hold for the $\kk G$-modules:
\begin{enumerate}
\item[(1)] $V_s\otimes V_{\chi}\cong V_{\chi}\otimes V_s\cong
V_{\tau^{-1}(s)}$.
\item[$(2)$] $V_s\otimes V_{\chi^{-t}}\cong V_{\tau^t(s)}$.
\item[$(3)$] There is a bijection
$\sigma_{s,t}$ from $V_s$ to $V_{\tau^t(s)}$ such that
$\sigma_{s,t}(hv)=\chi^t(h)h\sigma_{s,t}(v)$,
for any $h\in G$ and $v\in V_s$.
\item[$(4)$] $xV_s\cong V_{\tau(s)}$. Moreover, $V_i\cong V_j$ if and only if $xV_i\cong xV_j$, for $i,j\in\Omega$.
\item[$(5)$] The order of the permutation $\tau$ is $n$. Moreover, for any
    $s\in\Omega$, $x^t V_s\cong V_s$ if and only if $t$ is divisible by $n$.
\end{enumerate}
\end{lem}
\proof Part (1) and Part (2) are obvious. Part (3) is the same as \cite[Lemma 2.3]{WLZ}.

(4) The $\mathbbm{k}$-linear map $xv\mapsto u\otimes v$ gives an
isomorphism from $xV_s$ to $V_{\chi^{-1}}\otimes V_s$, where $0\neq
u\in V_{\chi^{-1}}$. Moreover, $V_i\cong V_j$ if and only if $V_{\tau(i)}\cong V_{\tau(j)}$ if and only if $xV_i\cong xV_j$, for
$i,j\in\Omega$.

(5) Note that the order of $\chi$ is $n$. It follows from Part (2) that the
order of $\tau$ is $n$ as well. Suppose the action of $g$ on $V_s$ is a scalar multiple by an element $\omega_s$. Then the action of $g$ on $x^t V_s$ is a scalar multiple by $\omega_sq^{-t}$.
If $x^tV_s\cong V_s$, then $\omega_s=\omega_sq^{-t}$, and hence $t$ is divisible by $n$ since the order of
$q$ is $n$. Conversely, if $t$ is divisible by $n$, it is obvious that $x^tV_s\cong V_s$ since the
order of $\tau$ is  $n$.\qed

Let $\langle\tau\rangle$ be the group generated by the permutation
$\tau$. Then $\langle\tau\rangle$ acts on the index sets
$\Omega_k$,$k\in\{0,1\}$, respectively. Let $\sim$ be the equivalence relation
on $\Omega_k$: $i\sim j$ if and only if $i$ and $j$ belong to the same $\langle\tau\rangle$-orbit. Denote by $[i]$ the equivalence class (or the $\langle\tau\rangle$-orbit) of $i$ and $\overline{\Omega}_k=\{[i]\mid i\in\Omega_k\}$ the set of all equivalence classes. By Lemma \ref{lemT3} (5), we have $|[i]|=n$ for any $i\in\Omega_k$, and hence $|\Omega_k|$ is divisible by $n$.

With the equivalence relation above, we have the following result.

\begin{prop}\label{propT11}
For any $i,j\in \Omega_1$, $P_i\cong P_j$ as $H$-modules (or equivalently, as $H(1-e)$-modules) if and only if $[i]=[j]$.
\end{prop}
\proof Suppose $P_i\cong P_j$ as $H$-modules, then the two are isomorphic as $\mathbbm{k}G$-modules. By
Krull-Schmidt theorem, the direct summand $V_i$ of $P_i$ is isomorphic to a direct summand $x^kV_j$ of $P_j$. It follows from Lemma \ref{lemT3} (4) that $[i]=[j]$. Conversely, if $[i]=[j]$, then $i=\tau^{k}(j)$ for some $k$. This implies that $V_i\cong V_{\tau^{k}(j)}\cong x^kV_j$ as $\mathbbm{k}G$-modules. Thus, $\dim P_i=\dim P_j$ since $\dim
V_i=\dim x^kV_j=\dim V_j$. Let $\varsigma$ denote the isomorphism from $V_i$ to $x^kV_j$. Then the $\mathbbm{k}$-linear map
$$P_i\rightarrow P_j,\ \sum_{s=0}^{n-1}x^sv_s\mapsto\sum_{s=0}^{n-1}y^s\varsigma(v_s),$$
 is an injective $H$-module morphism, where the term $y^s\varsigma(v_s)$ means the action of $y^s$ on $\varsigma(v_s)$. Since the dimensions of the two modules are equal,  we conclude that $P_i\cong P_j$.
\qed

Following Proposition \ref{propT11}, we obtain a partition on the set $\{P_j\mid j\in \Omega_1\}$. Let $P_{[j]}$ stand for a representative of the isomorphism class $[P_j]$ of $P_j$.
As a direct consequence of Theorem \ref{3.4} and Proposition \ref{propT11}, we obtain that
the set $\{P_{[j]}\mid j\in\Omega_1\}$ forms a complete set of simple $H(1-e)$-modules up to isomorphism.
We summarize the main result of this section as follows:

\begin{thm}\label{theorem}The set $\{M(k,i), P_{[j]}\mid i\in\Omega_0, 1\leq k\leq
n, j\in\Omega_1\}$ forms a complete set of finite dimensional indecomposable $H$-modules up to isomorphism.
\end{thm}

To end this section, we describe the dual of $H$-modules, which will be used later. For any finite dimensional $H$-module $M$, the dual $M^*:=\text{Hom}_{\kk}(M,\kk)$ is an $H$-module given by $(hf)(v)=f(S(h)v)$, for $h\in H,\ f\in M^*$ and $v\in M.$

\begin{prop}\label{3.7}
For any $i\in\Omega_0$, $j\in\Omega_1$ and $1\leq k\leq n$, the following hold:
\begin{enumerate}
  \item[(1)] $M(k,i)^*\cong M(k,\tau^{1-k}(i^*))$, where $i^*\in\Omega_0$ such that $(V_i)^*\cong V_{i^*}.$
  \item[(2)] $P_{[j]}^*\cong V_{\chi^{-1}}\otimes P_{[j^*]}$, where $j^*\in\Omega_1$ such that $(V_j)^*\cong V_{j^*}.$
  \item[(3)] $M(k,i)^{**}\cong M(k,i)$ and $P_{[j]}^{**}\cong P_{[j]}$.
\end{enumerate}
\end{prop}
\proof The proof of Part (2) is similar to the proof of Part (1), whereas the proof of Part (1) follows from \cite[Proposition 3.7 (1)]{WLZ}.
Part (3) is obvious because the square of the antipode of $H$ is inner.\qed

\section{Clebsch-Gordan formulas}

In this section, we compute the Clebsch-Gordan formulas for the
decompositions of tensor products of indecomposable
$H$-modules. For any two simple $\kk G$-modules $V_i$ and $V_j$, if $V_s$ is a direct summand of the tensor product $V_i\otimes V_j$, we let $\pi_s$  be the projection from $V_i\otimes V_j$ to $V_s$.

\begin{lem}\label{lemT4}
Let $V_s$ be a direct summand of $\mathbbm{k}G$-module $V_i\otimes V_j$.
\begin{enumerate}
\item If $i,j\in\Omega_0$, then $s\in\Omega_0$.
\item If $i\in\Omega_0$, $j\in\Omega_1$ or $i\in\Omega_1$, $j\in\Omega_0$, then
$s\in\Omega_1$.
\item If $i,j\in\Omega_1$, and $\lambda_i\lambda_j=1$, then $s\in\Omega_0$.
\item If $i,j\in\Omega_1$, and $\lambda_i\lambda_j\neq1$, then $s\in\Omega_1$.
\end{enumerate}
\end{lem}
\proof Follow the fact that $\lambda_i\lambda_j=\lambda_s$ because the projection $\pi_s:V_i\otimes V_j\rightarrow V_s$ exists. \qed

\begin{lem}\label{lemT5} Let $B_i$ be  a basis of $V_i$ for any $i\in\Omega$.
\begin{enumerate}
\item The set
$$\{y^s(x^tu\otimes v)\mid0\leq s,t\leq n-1,u\in B_i,v\in
B_j\}$$
forms a basis of $P_i\otimes P_j$, for any $i,j\in\Omega_1$.
\item The set $$\{y^s(x^tu\otimes v)\mid0\leq s\leq
n-1,0\leq t\leq k-1,u\in B_i,v\in B_j\}$$ forms a basis of $M(k,i)\otimes P_j$, for any $i\in\Omega_0$, $j\in\Omega_1$ and $1\leq k\leq n$.
\item The set
 $$\{y^s(v\otimes x^tu)\mid0\leq s\leq
n-1,0\leq t\leq k-1,u\in B_i,v\in B_j\}$$ forms a basis of
 $P_j\otimes M(k,i)$, for any $i\in\Omega_0$, $j\in\Omega_1$ and $1\leq k\leq n$.
\end{enumerate}
\end{lem}
\proof We only prove Part (1), and the proofs of  Part (2) and Part (3) are similar. Note that $\{x^tu\otimes x^sv\mid0\leq t,s\leq n-1,u\in
B_i,v\in B_j\}$ forms a basis of $P_i\otimes P_j$. It is sufficient to show that the following is true:
 \begin{equation}\label{nniltype}x^tu\otimes x^sv\in\textrm{sp}\{y^s(x^tu\otimes v)\mid0\leq s,t\leq n-1,u\in B_i,v\in
B_j\},\end{equation} for any $0\leq
s,t\leq n-1$, $u\in B_i$ and $v\in B_j$. We proceed by induction on $s$.
 It is obvious that (\ref{nniltype}) holds when $s=0$, for all $0\leq t\leq n-1$, $u\in B_i$ and $v\in B_j$. For a fixed $1\leq d\leq n-2$, suppose that (\ref{nniltype}) holds for $1\leq s\leq d$.
We consider the case $s=d+1$. Note that
$$\Delta(y^{d+1})=\sum_{p=0}^{d+1}\binom{d+1}{p}_q
y^{d+1-p}\otimes g^{d+1-p}y^p,$$
see e.g. \cite[Eq.(1)]{KR}. Then
\begin{align*}
&\ \ \ \ y^{d+1}(x^tu\otimes v)\\
&=\sum_{p=0}^{d+1}\binom{d+1}{p}_q
(y^{d+1-p}\otimes g^{d+1-p}y^p)(x^tu\otimes v)\\
&=\sum_{p=0}^{d}\binom{d+1}{p}_q
(y^{d+1-p}\otimes g^{d+1-p}y^p)(x^tu\otimes v)+(1\otimes y^{d+1})(x^tu\otimes v)\\
&=\sum_{p=0}^{d}\mu_p(x^{n_p}u\otimes x^pv)+(x^tu\otimes x^{d+1}v),
\end{align*}
where $\mu_p\in\mathbbm{k}$, and $n_p$ is the remainder after dividing $d+1-p+t$ by $n$.
By induction assumption, we have that the element $\sum_{p=0}^{d}\mu_p(x^{n_p}u\otimes x^pv)$ belongs to $\textrm{sp}\{y^s(x^tu\otimes v)\mid0\leq s,t\leq n-1,u\in B_i,v\in
B_j\}.$ This implies that
$$x^tu\otimes x^{d+1}v=y^{d+1}(x^tu\otimes v)-\sum_{p=0}^{d}\mu_p(x^{n_p}u\otimes x^pv)$$ belongs to
$\textrm{sp}\{y^s(x^tu\otimes v)\mid0\leq s,t\leq n-1,u\in B_i,v\in
B_j\}$, for any $0\leq t\leq n-1$, $u\in B_i$ and $v\in B_j$, as desired. \qed

With the bases of the tensor products of indecomposable modules given in Lemma \ref{lemT5}, we are able to describe the decomposition of the tensor product of any two indecomposable modules. We start with two indecomposable modules of form $P_i$, for $i\in \Omega_1$.

\begin{prop}\label{thT2}
For any $i,j\in\Omega_1$, we have the following decompositions:
\begin{enumerate}
\item[$(1)$] If $\lambda_i\lambda_j=1$, then $V_i\otimes V_j\cong\bigoplus_{s}V_s$ for some $s\in\Omega_0$. In this case, $$P_i\otimes
P_j\cong\bigoplus_{t=0}^{n-1}\bigoplus_{s}M(n,\tau^t(s)).$$
\item[$(2)$] If $\lambda_i\lambda_j\neq1$, then $V_i\otimes V_j\cong\bigoplus_{s}V_s$ for some $s\in\Omega_1$. In this case, $$P_i\otimes
P_j\cong\bigoplus_{t=0}^{n-1}\bigoplus_{s}P_{\tau^t(s)}.$$
\end{enumerate}
\end{prop}
\proof For any $0\leq t\leq n-1$, we let $\Gamma_t$ be the subset
$$\{x^tu\otimes v,y(x^tu\otimes
v),\cdots,y^{n-1}(x^tu\otimes v)\mid\textrm{for}\ u\in B_i,v\in
B_j\}$$
 of the basis $\{y^s(x^tu\otimes v)\mid0\leq s,t\leq n-1,u\in B_i,v\in
B_j\}$ of $P_i\otimes P_j$.
It is easy to see that the subspace $\textrm{sp}\Gamma_t$ is a
submodule of $P_i\otimes P_j$ with dimension $n\dim V_i\dim V_j$. Moreover, we have $P_i\otimes P_j\cong\bigoplus_{t=0}^{n-1}\textrm{sp}\Gamma_t$.

(1) If $\lambda_i\lambda_j=1$, then
$V_i\otimes V_j\cong\bigoplus_{s}V_s$, for some $s\in\Omega_0$ by Lemma \ref{lemT4}. For any $0\leq t\leq
n-1$, we consider the following $\mathbbm{k}$-linear map:
$$\varphi_t:\textrm{sp}\Gamma_t\rightarrow\bigoplus_{s}M(n,\tau^t(s)),\
y^k(x^tu\otimes v)\mapsto\sum_{s}x^k\sigma_{s,t}\pi_s(u\otimes
v),$$
for $0\leq k\leq n-1$, $u\in B_i$ and $v\in B_j,$ where the map $\sigma_{s,t}$ is given in Lemma \ref{lemT3} (3). It is tedious but straightforward to check that
the map $\varphi_t$ is an $H$-module isomorphism. It follows that
$$P_i\otimes P_j=\bigoplus_{t=0}^{n-1}\textrm{sp}\Gamma_t\cong\bigoplus_{t=0}^{n-1}\bigoplus_{s}M(n,\tau^t(s)).$$

(2) Suppose that  $\lambda_i\lambda_j\neq1$. By Lemma \ref{lemT4}, we have
$V_i\otimes V_j\cong\bigoplus_{s}V_s$, for some $s\in
\Omega_1$. Now for any $0\leq t\leq
n-1$, similar to Part (1), we have  the following $\mathbbm{k}$-linear map:
$$\psi_t:\textrm{sp}\Gamma_t\rightarrow\bigoplus_{s}P_{\tau^t(s)},\
y^k(x^tu\otimes
v)\mapsto\sum_{s}x^k\sigma_{s,t}\pi_s(u\otimes
v),$$ for $0\leq k\leq n-1$, $u\in B_i$ and $v\in B_j,$  which is in fact an $H$-module isomorphism.
Thus,
$$P_i\otimes P_j=\bigoplus_{t=0}^{n-1}\textrm{sp}\Gamma_t\cong\bigoplus_{t=0}^{n-1}\bigoplus_{s}P_{\tau^t(s)}.$$
We have completed the proof.\qed

Now we consider the decomposition of the tensor product of $M(k,i)$ with $P_j$.

\begin{prop}\label{thT3}
If $i\in\Omega_0$ and $j\in\Omega_1$, then $V_i\otimes V_j\cong\bigoplus_{s}V_s$ for some $s\in\Omega_1$. In this case, $$M(k,i)\otimes
P_j\cong P_j\otimes M(k,i)\cong\bigoplus_{t=0}^{k-1}\bigoplus_{s}P_{\tau^t(s)},$$ for any $1\leq k\leq n$.
\end{prop}
\proof Similar to the proof of Proposition \ref{thT2}.\qed

Note that the left $\overline{H}$-module category is a monoidal full subcategory of left $H$-module category by Proposition \ref{propT2}. The $H$-module decomposition of the tensor product $M(k,i)\otimes M(l,j)$ is the same as the $\overline{H}$-module decomposition. We present this decomposition directly as follows (see \cite[Proposition 4.2]{WLZ}):

\begin{prop}\label{propT7}
If $i,j\in\Omega_0$, then $V_i\otimes V_j\cong\bigoplus_{s}V_s$ for some $s\in\Omega_0$.
\begin{enumerate}
\item[$(1)$] If $k+l-1\leq n$, then
$$M(k,i)\otimes M(l,j)\cong\bigoplus_{s}\bigoplus_{t=0}^{\min\{k,l\}-1}M(k+l-1-2t,\tau^t(s)).$$
\item[$(2)$] If $k+l-1\geq n$, then
$$M(k,i)\otimes M(l,j)\cong\bigoplus_{s}(\bigoplus_{t=0}^{r}M(n,\tau^t(s))
\bigoplus\bigoplus_{t=r+1}^{\min\{k,l\}-1}M(k+l-1-2t,\tau^t(s))),$$
where $r=k+l-1-n$.
\end{enumerate}
\end{prop}

As an immediate consequence of Proposition \ref{thT2}, Proposition \ref{thT3} and Proposition \ref{propT7}, we have the following.
\begin{cor}\label{colT3}
We have that $M\otimes N\cong N\otimes M$, for any two finite dimensional $H$-modules $M$ and $N$.
\end{cor}

\section{The structure of Green rings}
Let $A$ be a Hopf algebra over a field $\kk$ and $F(A)$ the free abelian group generated by the isomorphism classes $[M]$ of finite dimensional
$A$-modules $M$. The abelian group $F(A)$ becomes a ring if we endow $F(A)$ with a multiplication given by the tensor product $[M][N]=[M\otimes N]$. The \emph{Green ring} (or \emph{representation ring}) $r(A)$ of the Hopf algebra $A$ is defined to be the quotient ring of $F(A)$ modulo the relations $[M\oplus N]=[M]+[N]$. The identity of the associative ring  $r(A)$ is represented by the trivial $A$-module $\kk$. Note that
$r(A)$ has a $\mathbb{Z}$-basis consisting of isomorphism classes of finite dimensional indecomposable $A$-modules (see e.g., \cite{Chen2, DH, LH, LZ}).

The \textit{Grothendieck ring}\index{Grothendieck ring} $G_0(A)$ of Hopf algebra $A$ is the quotient ring of $F(A)$ modulo short exact sequences of $A$-modules, i.e., $[Y]=[X]+[Z]$ if $0\rightarrow X\rightarrow Y\rightarrow Z\rightarrow0$ is exact. The Grothendieck ring  $G_0(A)$ possesses a basis given by isomorphism classes of simple $A$-modules.
Both $r(A)$ and  $G_0(A)$ are  augmented $\mathbb{Z}$-algebras with the dimension augmentation. Moreover, there is a natural ring epimorphism from $r(A)$ to $G_0(A)$ given by $$\phi:r(A)\rightarrow G_0(A),\ [M]\mapsto[M]=\sum_{[V]}[M:V][V],$$ where the sum runs over all finite dimensional simple $A$-modules and $[M:V]$ is the  multiplicity of $V$ in the composition series of $M$.

Let $r(\overline{H})$ and $r(H)$ be the Green rings of Hopf algebras
$\overline{H}$ and $H$ respectively. Then $r(H)$ is a commutative ring (Corollary \ref{colT3}) with the subring $r(\overline{H})$ (Proposition \ref{propT2}). Denote by $M[k,i]$ and $P_{[j]}$ the isomorphism classes of indecomposable modules $M(k,i)$ and $P_j$ respectively. We write $1$ for $[\kk]$ and $a$ for $[V_{\chi^{-1}}]$ respectively. Since the order of $\chi$ is $n$, we have $a^n=1$.

\begin{prop}\label{propT3}
For any $i,j\in\Omega_1$, the following hold in $r(H)$:
\begin{enumerate}
\item[$(1)$] If $\lambda_i\lambda_j=1$, then $V_i\otimes V_j\cong\bigoplus_{s}V_s$ for some $s\in\Omega_0$. In this case,
$$P_{[i]}P_{[j]}=(1+a+\cdots+a^{n-1})\sum_{s}M[n,s].$$
\item[$(2)$] If $\lambda_i\lambda_j\neq1$, then $V_i\otimes V_j\cong\bigoplus_{s}V_s$ for some $s\in\Omega_1$. In this case,
$$P_{[i]}P_{[j]}=n\sum_{s}P_{[s]}.$$
\end{enumerate}
\end{prop}
\proof (1) It follows from Proposition \ref{thT2} (1) that
\begin{align*}
P_{[i]}P_{[j]}&=\sum_{t=0}^{n-1}\sum_{s}M[n,\tau^t(s)]=\sum_{t=0}^{n-1}\sum_{s}a^{t}M[n,s]\\
&=(1+a+\cdots+a^{n-1})\sum_{s}M[n,s].
\end{align*}
In addition, the foregoing expression  is well-defined. Indeed, if $i_1=\tau^k(i)$ and $j_1=\tau^p(j)$, for some
integers $k$ and $p$, then $V_{i_1}\cong V_i\otimes V_{\chi^{-k}}$ and $V_{j_1}\cong V_j\otimes V_{\chi^{-p}}$. These yield that $\lambda_{i_1}\lambda_{j_1}=\lambda_{i}\chi^{-k}(g^n)\lambda_{j}\chi^{-p}(g^n)=\lambda_i\lambda_j=1$ and
\begin{align*}
V_{i_1}\otimes V_{j_1}
&\cong V_i\otimes V_{\chi^{-k}}\otimes V_{j}\otimes V_{\chi^{-p}}\cong V_i\otimes V_j\otimes V_{\chi^{-(k+p)}}\\
&\cong(\bigoplus_{s}V_s)\otimes
V_{\chi^{-(k+p)}}\cong\bigoplus_{s}V_{\tau^{k+p}(s)}.
\end{align*} Note that $a^n=1$. By Proposition \ref{thT2} (1), we have
\begin{align*}
P_{[i_1]}P_{[j_1]}&=\sum_{t=0}^{n-1}\sum_{s}M[n,\tau^t(\tau^{k+p}(s))]=\sum_{t=0}^{n-1}\sum_{s}a^{t+k+p}M[n,s]\\
&=(1+a+\cdots+a^{n-1})\sum_{s}M[n,s]=P_{[i]}P_{[j]}.
\end{align*}

(2) Note that $P_{[{\tau^t(s)}]}=P_{[s]}$ since $\tau^t(s)\sim
s$. By Proposition \ref{thT2} (2), we have
$$P_{[i]}P_{[j]}=\sum_{t=0}^{n-1}\sum_{s}P_{[{\tau^t(s)}]}
=\sum_{t=0}^{n-1}\sum_{s}P_{[s]}=n\sum_{s}P_{[s]}.$$
To show that this expression is well-defined, we assume that $i_1=\tau^k(i)$ and
$j_1=\tau^p(j)$, for some integers $k$ and $p$. Then $\lambda_{i_1}\lambda_{j_1}\neq1$ since $\lambda_i\lambda_j\neq1$, and
$V_{i_1}\otimes V_{j_1}\cong\bigoplus_{s}V_{\tau^{k+p}(s)}$ as shown
above. It follows from Proposition \ref{thT2} (2) that
$$P_{[i_1]}P_{[j_1]}=\sum_{t=0}^{n-1}\sum_{s}P_{[{\tau^t(\tau^{k+p}(s))}]}
=\sum_{t=0}^{n-1}\sum_{s}P_{[s]}=n\sum_{s}P_{[s]}=P_{[i]}P_{[j]}.$$
The proof is completed. \qed

\begin{prop}\label{propT4}
If $i\in\Omega_0$ and $j\in\Omega_1$, then $V_i\otimes V_j\cong\bigoplus_{s}V_s$ for some $s\in\Omega_1$. In this case,
$$M[k,i]P_{[j]}=k\sum_{s}P_{[s]}$$ for any $1\leq k\leq n$.
In particular, $[V_i]P_{[j]}=\sum_{s}P_{[s]}$,
$M[2,0]P_{[j]}=2P_{[j]}$ and $(1+a-M[2,0])P_{[j]}=0$.
\end{prop}
\proof By Proposition \ref{thT3}, we have
$$M[k,i]P_{[j]}=\sum_{t=0}^{k-1}\sum_{s}P_{[\tau^t(s)]}
=\sum_{t=0}^{k-1}\sum_{s}P_{[s]}
=k\sum_{s}P_{[s]}.$$ The proof that the above expression is well-defined is similar to those in Proposition \ref{propT3} .
\qed


In view of the relations given in Proposition \ref{propT3} and Proposition \ref{propT4}, it is difficult to present the Green ring $r(H)$ of $H$ in terms of (minimal) generators and relations as done in \cite{Chen, HOYZ, LZ, Wa}. However, we are still able to describe the structure of $r(H)$ in terms of the Green ring $r(\overline{H})$.

Let $\mathcal {P}$ be the free abelian group generated by isomorphism classes of finite dimensional indecomposable projective $H$-modules. That is,
$$\mathcal {P}=\mathbb{Z}\{M[n,i],P_{[j]}\mid i\in\Omega_0,\ j\in\Omega_1\}.$$
Then $\mathcal {P}$ is an ideal of $r(H)$ since any $H$-module tensors with a projective module is also projective.
The direct sum $r(\overline{H})\bigoplus\mathcal {P}$ is a commutative ring with the multiplication given by $$(b_1,c_1)(b_2,c_2)=(b_1b_2,b_1c_2+c_1b_2+c_1c_2),$$ for any $b_1,b_2\in r(\overline{H})$ and $c_1,c_2\in\mathcal {P}$. Obviously, $r(\overline{H})\bigoplus\mathcal {P}$ can be regarded as a certain trivial extension of $r(\overline{H})$ with respect to $\mathcal {P}$.

\begin{thm}\label{th3.6.3}
Let $\mathcal {I}$ be the submodule of $\mathbb{Z}$-module $r(\overline{H})\bigoplus\mathcal {P}$ generated by the elements $(-M[n,i],M[n,i])$ for $i\in\Omega_0$. Then $\mathcal {I}$ is an ideal of $r(\overline{H})\bigoplus\mathcal {P}$ and the quotient ring $(r(\overline{H})\bigoplus\mathcal {P})/\mathcal {I}$ is isomorphic to $r(H)$.
\end{thm}
\proof For any $b\in r(\overline{H})$ and $c\in\mathcal {P}$, we have $$(b,c)(-M[n,i],M[n,i])=(-bM[n,i],bM[n,i])\in\mathcal {I}$$
since $bM[n,i]$ is always a $\mathbb{Z}$-linear combination of elements of the form $M[n,j]$ for $j\in\Omega_0$. Thus, $\mathcal {I}$ is a two-sided ideal of $r(\overline{H})\bigoplus\mathcal {P}$ since $r(\overline{H})\bigoplus\mathcal {P}$ is commutative. Note that the $\mathbb{Z}$-linear map from $r(\overline{H})\bigoplus\mathcal {P}$ to $r(H)$ given by $(b,c)\mapsto b+c$ is a ring epimorphism with the kernel equal to $\mathcal {I}$. We conclude that $(r(\overline{H})\bigoplus\mathcal {P})/\mathcal {I}\cong r(H)$, as desired.
\qed

\section{Grothendieck rings and Cartan matrices}

In this section, we study the Grothendieck ring $G_0(H)$ of $H$ in terms of the Green ring $r(H)$. We then describe the Cartan map from $\mathcal {P}$ to the Grothendieck ring $G_0(H)$ and obtain the Cartan matrix with respect to given bases.

Let $M$ be a finite dimensional $H$-module. Recall that the $\kk$-linear dual $M^*$ is an $H$-module with the $H$-module structure given by Proposition \ref{3.7}. This leads to an anti-automorphism $*$ of the Green ring $r(H)$: $[M]^*:=[M^*]$. Since $M^{**}\cong M$, the $*$-operator is an involution of $r(H)$.

Now for any indecomposable $H$-module $Z$, if $Z$ is not projective, there exists a unique almost split sequence $0\rightarrow X\rightarrow Y\rightarrow Z\rightarrow0$ ending at $Z$. We take the  same notation introduced  in \cite[Section 4, ChVI]{ARS} (see also \cite{WLZ}) and denote by $\delta_{[Z]}$ the element $[X]-[Y]+[Z]$ in $r(H)$.  If $Z$ is projective, then we write $\delta_{[Z]}$ for $[Z]-[\textrm{rad}Z].$

For any two indecomposable $H$-modules $X$ and $Y$, similar to the one in \cite{WLZ}, we define:
$$([X],[Y])=\dim\text{Hom}_H(X,Y^*),$$
which extends to an associative, symmetric and  non-degenerate $\mathbb{Z}$-bilinear form on $r(H)$. This bilinear form  satisfies the following:
\begin{equation}\label{3.8}(\delta^*_{[X]},[Y])=1\ \text{if}\ X\cong Y,\ \text{and}\ 0\ \text{otherwise}.\end{equation}
It follows that $r(H)$ is a symmetric algebra over $\mathbb{Z}$ with a pair of dual bases
$$\{M[k,i],\ P_{[j]}\mid i\in\Omega_0,\ j\in\Omega_1,\ 1\leq k\leq n\}$$ and $$\{\delta^*_{M[k,i]},\ P_{[j]}^*\mid i\in\Omega_0,\ j\in\Omega_1,\ 1\leq k\leq n\}.$$

With the above pair of dual bases, we can express any element $x$ in $r(H)$  as follows:
\begin{equation}\label{equ4.8}
x=\sum_{i\in\Omega_0}\sum_{k=1}^n(\delta_{M[k,i]}^*,x)M[k,i]+\sum_{[j]\in\overline{\Omega}_1}(P_{[j]}^*,x)P_{[j]},
\end{equation}
or equivalently,
\begin{equation}\label{equ4.9}
x=\sum_{i\in\Omega_0}\sum_{k=1}^n(x,M[k,i])\delta_{M[k,i]}^*+\sum_{[j]\in\overline{\Omega}_1}(x,P_{[j]})P_{[j]}^*.
\end{equation}

\begin{lem}\label{lemT1}
The following equations hold in $r(H)$:
\begin{enumerate}
  \item[(1)] $\delta_{M[k,i]}=\delta_{[\kk]}M[k,i]$ and $\delta_{[P_j]}=P_{[j]}$, where $\delta_{[\kk]}=1+a-M[2,0]$, $i\in\Omega_0$, $j\in\Omega_1$ and $1\leq k\leq n-1$.
  \item[(2)] $\delta_{[\kk]}M[n,i]=\delta_{[\kk]}P_{[j]}=0$, for $i\in\Omega_0$ and $j\in\Omega_1$.
\end{enumerate}
\end{lem}
\proof (1) Note that an almost split sequence of $\overline{H}$-modules is also an almost split sequence of $H$-modules since the $\overline{H}$-module category is a monoidal full subcategory of the $H$-module category. It follows that the almost split sequence of $\overline{H}$-modules (see \cite[Proposition 3.5]{WLZ})
$$0\rightarrow V_{\chi^{-1}}\otimes M(k,i)\rightarrow M(2,0)\otimes M(k,i)\rightarrow M(k,i)\rightarrow0$$
is almost split as $H$-modules, for any $i\in\Omega_0$ and $1\leq k\leq n-1$. In particular, the almost split sequence ending at $M(1,0)=\kk$ is
$$0\rightarrow V_{\chi^{-1}}\rightarrow M(2,0)\rightarrow \kk\rightarrow0.$$
Consequently, we have $\delta_{[\kk]}=1+a-M[2,0]$ and $\delta_{M[k,i]}=\delta_{[\kk]}M[k,i],$ for $i\in\Omega_0$ and $1\leq k\leq n-1$.
It is obvious that $\delta_{[P_j]}=[P_{j}]-[\text{rad}P_j]=P_{[j]}$ since $P_{j}$ is simple projective for any $j\in\Omega_1$.

(2) Follows from the fact that an almost split sequence tensoring  with a projective module is split.
\qed

The natural ring epimorphism from the Green ring $r(H)$ to the Grothendieck ring $G_0(H)$ of $H$ is as follows:
\begin{equation*}\phi:r(H)\rightarrow G_0(H),\ M[k,i]\mapsto[V_i](1+a+\cdots+a^{k-1}),\ P_{[j]}\mapsto P_{[j]}\end{equation*} for $i\in\Omega_0,$ $j\in\Omega_1$ and $1\leq k\leq n$. The kernel of $\phi$ is exactly the free abelian group generated by almost split sequences ending at indecomposable non-projective modules (cf. \cite[Theorem 4.4, ChVI]{ARS}). That is,
$$\ker\phi=\mathbb{Z}\{\delta_{M[k,i]}\mid i\in\Omega_0,\ 1\leq k\leq n-1\}.$$
By Lemma \ref{lemT1}, $\ker\phi$ is a principal ideal of $r(H)$ generated by $\delta_{[\kk]}$. Moreover, $(\ker\phi)^*=\ker\phi$ is $*$-invariant since $\delta_{[\kk]}^*=a^{-1}\delta_{[\kk]}$.

Thanks to the map $\phi$ given above, one is able to describe the multiplication rule in $G_0(H)$.
Since $\{[V_i],P_{[j]}\mid i\in\Omega_0,j\in\Omega_1\}$ is a basis of $G_0(H)$, it is enough to look at the products of basis elements. The products $[V_i][V_j]$ and $[V_i]P_{[k]}$ in $G_0(H)$ are the same as in $r(H)$. However, the product of $P_{[i]}$ and $P_{[j]}$ in $G_0(H)$ becomes
$$P_{[i]}P_{[j]}=\begin{cases}
n(1+a+\cdots+a^{n-1})\sum_s[V_s], & \lambda_i\lambda_j=1,\\
n\sum_sP_{[s]}, & \lambda_i\lambda_j\neq1,
\end{cases}
$$ where the index $s$ is determined by the decomposition $V_i\otimes V_j\cong\bigoplus_s V_s$ as $\kk G$-modules. We refer to Proposition \ref{propT3} for their multiplication in $r(H)$.

Denote by $$\mathcal {P}^{\perp}=\{x\in r(H)\mid(x,y)=0,\ \forall y\in\mathcal {P}\},$$ the orthogonal ideal of $\mathcal{P}$ with respect to the form $(-,-)$. We are ready to describe the structure of the Grothendieck ring $G_0(H)$ as follows.

\begin{prop}\label{prop4.50}Let $r(\kk G)$ be the Green ring of the group algebra $\kk G$.
\begin{enumerate}
  \item We have $\ker\phi=\mathcal {P}^{\perp}$, hence $G_0(H)\cong r(H)/\mathcal {P}^{\perp}$.
  \item The Grothendieck ring $G_0(H)$ is isomorphic to the subring of $r(\kk G)$ generated by $[V_i]$ and $(1+a+\cdots+a^{n-1})[V_j]$, for any $i\in\Omega_0$ and $j\in\Omega_1$. 
\end{enumerate}
\end{prop}
\proof (1) By Lemma \ref{lemT1} (2), $(\delta_{[\kk]},P_{[j]})=([\kk],\delta_{[\kk]}P_{[j]})=([\kk],0)=0$. Similarly, we have $(\delta_{[\kk]},M[n,i])=0$.
These imply that $\delta_{[\kk]}\in\mathcal {P}^{\perp}$, and hence $\ker\phi\subseteq\mathcal {P}^{\perp}$.
Conversely, for any $x\in\mathcal {P}^{\perp}$, by (\ref{equ4.9}), \begin{align*}x&=\sum_{i\in\Omega_0}\sum_{k=1}^n(x,M[k,i])\delta_{M[k,i]}^*+\sum_{[j]\in\overline{\Omega}_1}(x,P_{[j]})P_{[j]}^*\\
&=\sum_{i\in\Omega_0}\sum_{k=1}^{n-1}(x,M[k,i])\delta_{M[k,i]}^{*}\ \ (\text{as}\ x\in\mathcal {P}^{\perp}).
\end{align*}
It follows that $x^*=\sum_{i\in\Omega_0}\sum_{k=1}^{n-1}(x,M[k,i])\delta_{M[k,i]}\in\ker\phi$. Hence, $x=x^{**}\in(\ker\phi)^*=\ker\phi$, and $\mathcal {P}^{\perp}\subseteq\ker\phi$.

(2) Consider the $\mathbb{Z}$-linear map $\varphi$ from $G_0(H)$ to $r(\mathbbm{k}G)$ given by
$$\varphi ([V_i])=[V_i]\ \text{and}\ \varphi(P_{[j]})=(1+a+\cdots+a^{n-1})[V_j],$$
for any $i\in\Omega_0$ and $j\in\Omega_1$. The map $\varphi$ is well-defined since the sum $(1+a+\cdots+a^{n-1})[V_j]$ is a $\langle\tau\rangle$-orbit sum and $(1+a+\cdots+a^{n-1})[V_j]=(1+a+\cdots+a^{n-1})[V_k]$ if and only if $[j]=[k]$.
It is straightforward to verify that the map $\varphi$ is a ring homomorphism. We claim that the map is injective. Indeed, if $\sum_{i\in\Omega_0}\alpha_i[V_i]+\sum_{[j]\in\overline{\Omega}_1}\beta_{[j]}P_{[j]}$ belongs to the kernel of $\varphi$, then
$\sum_{i\in\Omega_0}\alpha_i[V_i]+\sum_{[j]\in\overline{\Omega}_1}\beta_{[j]}(1+a+\cdots+a^{n-1})[V_j]=0$. This yields that each $\alpha_i=0$ and $\beta_{[j]}=0$ since the elements of the set $\{[V_i],(1+a+\cdots+a^{n-1})[V_j]\mid i\in\Omega_0,[j]\in\overline{\Omega}_1\}$ are $\mathbb{Z}$-linear independent.
\qed

\begin{rem}\label{rem4.51}
It was shown in \cite{Lo} that if the Jacobson radical of a finite dimensional Hopf algebra is a Hopf ideal, then its Grothendieck ring is semiprime. The Proposition \ref{prop4.50} (2) means that the Grothendieck ring $G_0(H)$ of $H$ is semiprime although the Jacobson radical of $H$ is not a Hopf ideal.
\end{rem}

The map $\phi$ restricting to the the ideal $\mathcal {P}$ of $r(H)$ gives rise to the Cartan map as follows:
\begin{equation*}\phi|_{\mathcal {P}}:\mathcal {P}\rightarrow G_0(H),\ M[n,i]\mapsto[V_i](1+a+\cdots+a^{n-1}),\ P_{[j]}\mapsto P_{[j]}\end{equation*} for $i\in\Omega_0$ and $j\in\Omega_1$. Note that, for any $i\in\Omega_0$, the $\langle\tau\rangle$-orbit containing $i$ is $[i]=\{i,\tau(i),\cdots,\tau^{n-1}(i)\}$. Accordingly, we write  $V_{[i]}$ for the set $\{[V_i],[V_{\tau(i)}],\cdots,[V_{\tau^{n-1}(i)}]\}$. The disjoint union
$$(\cup_{[i]\in\overline{\Omega}_0}V_{[i]})\cup\{P_{[j]}\mid[j]\in\overline{\Omega}_1\}$$ forms a basis of $G_0(H)$ and the set of their projective covers $$(\cup_{[i]\in\overline{\Omega}_0}P_{V_{[i]}})\cup\{P_{[j]}\mid[j]\in\overline{\Omega}_1\}$$ forms a basis of $\mathcal {P}$, where $P_{V_{[i]}}=\{M[n,i],M[n,\tau(i)],\cdots,M[n,\tau^{n-1}(i)]\}$. Note that the map $\phi|_{\mathcal {P}}$ maps any element in $P_{V_{[i]}}$ into the same element $[V_i](1+a+\cdots+a^{n-1})$. It follows that the Cartan matrix of the map $\phi|_{\mathcal {P}}$ with respect to the above bases is of a block-diagonal form with $|\overline{\Omega}_0|+1$ blocks:
$$\left(
    \begin{array}{cccc}
      \mathbf{E} &  &  &  \\
       & \ddots &  &  \\
       &  & \mathbf{E} &  \\
       &  &  & \bf{I} \\
    \end{array}
  \right),
$$
where $\mathbf{E}$ is a square matrix of order $n$ with each entry 1 and $\mathbf{I}$ is the identity matrix of order $|\overline{\Omega}_1|$.

\section{Jacobson radicals and idempotents of Green rings}
In this section, we shall use the Cartan map of $H$ to describe the Jacobson radical and the idempotents of $r(H)$. We show that the  Jacobson radical of $r(H)$ is exactly the kernel of the Cartan map generated by one element, and $r(H)$ has no non-trivial idempotents. We first need the following lemma.

\begin{lem}\label{19.2} The following hold in $r(H)$:
\begin{enumerate}
  \item[(1)] For any two indecomposable $H$-modules $X$ and $Y$, the coefficient of  $[\kk]$ in the expression of $[Y][X]^*$ is 1 if $X\cong Y$ and neither of them is projective. Otherwise, the coefficient is 0.
  \item[(2)] If $xx^*=0$ for $x\in r(H)$, then $x\in\mathcal {P}$.
\end{enumerate}
\end{lem}
\proof (1) If $X$ is projective, then the coefficient of $[\kk]$ in $[Y][X]^*$ is 0 since $\kk$ itself is not projective. If $X$ is not projective, by (\ref{3.8})  the coefficient of $[\kk]$ in $[Y][X]^*$ is $(\delta^*_{[\kk]},[Y][X]^*)=(\delta^*_{[X]},[Y])$, which is equal to 1 if $X\cong Y$, and to 0 otherwise.

(2) Suppose $x=\sum_{i\in\Omega_0}\sum_{1\leq k\leq n-1}\alpha_{ki}M[k,i]+x_0,$
for $\alpha_{ki}\in\mathbb{Z}$ and $x_0\in\mathcal {P}$. By Part (1), the coefficient of $[\kk]$ in $xx^*$ is $\sum_{i\in\Omega_0}\sum_{1\leq k\leq n-1}\alpha_{ki}^2$. Thus, if $xx^*=0$, then $\alpha_{ki}=0$ for any $i\in\Omega_0$ and $1\leq k\leq n-1$. This implies that $x=x_0\in\mathcal {P}$.
\qed

Note that the Green ring $r(H)$ is commutative and finitely generated as an algebra over $\mathbb{Z}$. This means that $r(H)$ is a commutative Jacobson ring. In this case, the Jacobson radical of $r(H)$ is the set consisting of all nilpotent elements of the ring.
\begin{thm}\label{th21}
The Jacobson radical of $r(H)$ is $J(r(H))=\mathcal {P}\cap\mathcal {P}^{\perp}=\ker(\phi|_{\mathcal {P}})$.
\end{thm}
\proof The equality $\mathcal {P}\cap\mathcal {P}^{\perp}=\ker(\phi|_{\mathcal {P}})$ is obvious since $\mathcal {P}^{\perp}=\ker\phi$ by Proposition \ref{prop4.50}. In the following, we verify that $J(r(H))=\mathcal {P}\cap\mathcal {P}^{\perp}$.

For any $x\in\mathcal {P}\cap\mathcal {P}^{\perp}$ and $y\in r(H)$, we have $xy\in\mathcal {P}\cap\mathcal {P}^{\perp}$, implying that $(x^2,y)=(x,xy)=0$. This means that $x^2=0$ since the bilinear form $(-,-)$ on $r(H)$ is non-degenerate. Consequently, $x\in J(r(H))$ and $\mathcal {P}\cap\mathcal {P}^{\perp}\subseteq J(r(H))$.

Conversely, for any $x\in J(r(H))$, we have $x\in\ker\phi$ since $\phi(x)$ is a nilpotent element of $G_0(H)$, whereas $G_0(H)$ has no non-trivial nilpotent elements since it is semiprime (Remark \ref{rem4.51}). On the other hand, we claim that $x\in\mathcal {P}$. In fact, we let $$x=\sum_{i\in\Omega_0}\sum_{1\leq j\leq n-1}\alpha_{ji}M[j,i]+x_0,$$ for $\alpha_{ji}\in\mathbb{Z}$ and $x_0\in\mathcal {P}$. Then $$xx^*=\sum_{i,k\in\Omega_0}\sum_{1\leq j,l\leq n-1}\alpha_{ji}\alpha_{lk}M[j,i]M[l,k]^*+x_1$$ for some $x_1\in\mathcal {P}$ since $\mathcal {P}$ is an ideal of $r(H)$ satisfying $\mathcal {P}=\mathcal {P}^*$. Denote by
$$y:=xx^*=\sum_{i\in\Omega_0}\sum_{1\leq j\leq n-1}\beta_{ji}M[j,i]+x_2,$$ for some $x_2\in\mathcal {P}$. By Lemma \ref{19.2} (1), the coefficient of $[\kk]$ in $y$ is $\beta_{10}=\sum_{i\in\Omega_0}\sum_{1\leq j\leq n-1}\alpha_{ji}^2$.
Consider $$y^2=yy^*=\sum_{i,k\in\Omega_0}\sum_{1\leq j,l\leq n-1}\beta_{ji}\beta_{lk}M[j,i]M[l,k]^*+x_3$$ for some $x_3\in\mathcal {P}$.
Again by Lemma \ref{19.2} (1), the coefficient of $[\kk]$ in $y^2$ is $\sum_{i\in\Omega_0}\sum_{1\leq j\leq n-1}\beta_{ji}^2$.
If $\beta_{10}\neq0$, then $\sum_{i\in\Omega_0}\sum_{1\leq j\leq n-1}\beta_{ji}^2\neq0$, and hence $y^2\neq0$. Repeating this process, we obtain that if $\beta_{10}\neq0$, then $y^{2^n}\neq0$ for any $n>0$. This contradicts to the fact that $y\in J(r(H))$. Thus, $\beta_{10}=0$ and  $x=x_0\in\mathcal {P}$. We obtain that $J(r(H))\subseteq \ker\phi\cap\mathcal {P}=\mathcal {P}\cap\mathcal {P}^{\perp}$. The proof is completed. \qed

In the following, we use Theorem \ref{th21} to describe the Jacobson radical of $r(H)$ in terms of generators.

\begin{thm}\label{th4.3.4}
The Jacobson radical of $r(H)$ is a principal ideal generated by $(1-a)M[n,0]$.
\end{thm}
\proof Note that $(1-a)M[n,0]\in\ker(\phi|_{\mathcal {P}})$. The ideal of $r(H)$ generated by $(1-a)M[n,0]$ is contained in $\ker(\phi|_{\mathcal {P}})$. Conversely, for any $\sum_{i\in\Omega_0}\alpha_iM[n,i]+\sum_{[j]\in\overline{\Omega}_1}\beta_{[j]}P_{[j]}$ in $\ker(\phi|_{\mathcal {P}})$, we have
$$\sum_{i\in\Omega_0}\alpha_i(1+a+\cdots+a^{n-1})[V_i]+\sum_{[j]\in\overline{\Omega}_1}\beta_{[j]}P_{[j]}=0.$$ It follows that $\beta_{[j]}=0$ for any $[j]\in\overline{\Omega}_1$, and  \begin{equation}\label{equaaa}\sum_{i\in\Omega_0}\alpha_i(1+a+\cdots+a^{n-1})[V_i]=0.\end{equation} Note that $(1+a+\cdots+a^{n-1})[V_i]=(1+a+\cdots+a^{n-1})[V_j]$ if and only if $[i]=[j]$. Thus, the equality (\ref{equaaa}) can be written as $$\sum_{[i]\in\overline{\Omega}_0}(\alpha_i+\alpha_{\tau(i)}+\cdots+\alpha_{\tau^{n-1}(i)})(1+a+\cdots+a^{n-1})[V_i]=0.$$ This implies that $\alpha_i+\alpha_{\tau(i)}+\cdots+\alpha_{\tau^{n-1}(i)}=0$ since the $\langle\tau\rangle$-orbit sum $(1+a+\cdots+a^{n-1})[V_i]$ for $[i]\in\overline{\Omega}_0$ are $\mathbb{Z}$-linear independent. So far we have verified that any element of $\ker(\phi|_{\mathcal {P}})$ is of the form $\sum_{i\in\Omega_0}\alpha_iM[n,i],$ where the coefficients satisfy $\alpha_i+\alpha_{\tau(i)}+\cdots+\alpha_{\tau^{n-1}(i)}=0$, for $i\in\Omega_0$. The equality $\alpha_i+\alpha_{\tau(i)}+\cdots+\alpha_{\tau^{n-1}(i)}=0$ implies the following equality:
$$\alpha_i+\alpha_{\tau(i)}a+\cdots+\alpha_{\tau^{n-1}(i)}a^{n-1}=\alpha_{\tau(i)}(a-1)+\cdots+\alpha_{\tau^{n-1}(i)}(a^{n-1}-1),$$ where the right hand side is in the ideal of $r(H)$ generated by $1-a$. As a result, we obtain that the element
\begin{align*}
\sum_{i\in\Omega_0}\alpha_iM[n,i]&=(\sum_{i\in\Omega_0}\alpha_i[V_i])M[n,0]\\
&=\sum_{[i]\in\overline{\Omega}_0}(\alpha_i[V_i]+\alpha_{\tau(i)}[V_{\tau(i)}]+\cdots+\alpha_{\tau^{n-1}(i)}[V_{\tau^{n-1}(i)}])M[n,0]\\
&=\sum_{[i]\in\overline{\Omega}_0}(\alpha_i+\alpha_{\tau(i)}a+\cdots+\alpha_{\tau^{n-1}(i)}a^{n-1})[V_i]M[n,0],
\end{align*}
sits in the ideal of $r(H)$ generated by $(1-a)M[n,0]$. The proof is completed.
\qed

We have shown that the Jacobson radical $J(r(H))$ of $r(H)$ is precisely $\mathcal {P}\cap\mathcal {P}^{\perp}$. So the generator $M[n,0](1-a)$ of $J(r(H))$ belongs to $\mathcal {P}\cap\mathcal {P}^{\perp}$. In particular, $M[n,0](1-a)$ is in $\mathcal {P}^{\perp}$, a principal ideal of $r(H)$ generated by $\delta_{[\kk]}$. In other words, the element $M[n,0](1-a)$ can be written as the product of $\delta_{[\kk]}$ with some element of $r(H)$. Indeed, by induction on $k$, we have $\delta_{[\kk]}(M[1,0]+\cdots+M[k,0])=1+aM[k,0]-M[k+1,0],$ for $1\leq k\leq n-1$. This implies that
\begin{align*}&\ \ \ \ \delta_{[\kk]}\sum_{k=1}^{n-1}(M[1,0]+\cdots+M[k,0])a^{n-1-k}\\
&=\sum_{k=1}^{n-1}(1+aM[k,0]-M[k+1,0])a^{n-1-k}\\
&=(1+a+\cdots+a^{n-1})-M[n,0].
\end{align*}
Thus, $M[n,0]$ can be written as
\begin{equation*}M[n,0]=(1+a+\cdots+a^{n-1})-\delta_{[\kk]}\sum_{k=1}^{n-1}(M[1,0]+\cdots+M[k,0])a^{n-1-k}.\end{equation*}
Multiplying both sides of the equation by $(1-a)$, we obtain:
$$M[n,0](1-a)=-\delta_{[\kk]}\sum_{k=1}^{n-1}(M[1,0]+\cdots+M[k,0])a^{n-1-k}(1-a).$$

To end this section, we show that the Green ring $r(H)$ has no non-trivial idempotents. First, observe that the Green ring of a group algebra possesses such a property. This fact might be found in other literature which we don't have at hands. So we include it in the following remark.

\begin{rem}\label{prop4.3}
For any finite group $G$, the Green ring (i.e., Grothendieck ring) $r(\kk G)$ of the group algebra $\kk G$ has no non-trivial idempotents.
Indeed, suppose that $\{V_i\mid i\in\Omega\}$ is a complete set of simple $\kk G$-modules up to isomorphism. For any $x\in r(\kk G)$, if we write $x=\sum_{i\in\Omega}\alpha_i[V_i]$ for $\alpha_i\in\mathbb{Z}$, then the coefficient of $[\kk]$ in $xx^*$ is $\sum_{i\in\Omega}\alpha_i^2$. Thus, $x=0$ if and only if $xx^*=0$. If $E$ is a primitive idempotent of $r(\kk G)$, so is $E^*$. Since $r(\kk G)$ is commutative and the duality operator $*$ is an anti-automorphism of $r(\kk G)$, either $E=E^*$ or $EE^*=0$. If $EE^*=0$, it follows that $E=0$. If $E=E^*$, comparing the coefficient of $[\kk]$ in both sides of the equation $EE^*=E$, we obtain that $E=0$ or $E=1$.
\end{rem}

\begin{thm}\label{th1.3}
The Green ring $r(H)$ has no non-trivial idempotents.
\end{thm}
\proof Let $E$ be a primitive idempotent of $r(H)$. We first prove that $E\in\mathcal {P}$
or $1-E\in\mathcal {P}$. Note that $E$ is primitive. So is $E^*$. Then $E=E^*$ or $EE^*=0$. If $EE^*=0$, by Lemma \ref{19.2} (2), $E\in\mathcal {P}$. If $E=E^*$, let $E=\sum_{i\in\Omega_0}\sum_{1\leq k\leq n-1}\alpha_{ki}M[k,i]+E_0,$ for $\alpha_{ki}\in\mathbb{Z}$ and $E_0\in\mathcal {P}$. Comparing the coefficient of $[\kk]$ in both sides of the equation
$EE^*=E$, we obtain that $E=\alpha_{10}[\kk]+E_0$, where $\alpha_{10}=0$ or 1. Therefore $E\in\mathcal {P}$ or $1-E\in\mathcal {P}$.
If $E\in\mathcal {P}$, let $E=\sum_{i\in\Omega_0}\alpha_{ni}M[n,i]+\sum_{[j]\in\overline{\Omega}_1}\alpha_{[j]}P_{[j]}$. Then
$$\phi(E)=(\sum_{i\in\Omega_0}\alpha_{ni}[V_i])(1+a+\cdots+a^{n-1})+\sum_{[j]\in\overline{\Omega}_1}\alpha_{[j]}P_{[j]}.$$
Note that $\phi(E)$ is an idempotent of $G_0(H)\subseteq r(\kk G)$ since $E$ is an idempotent of $r(H)$.
By Remark \ref{prop4.3},  $\phi(E)$ is equal to $0$ or $1$. However, $\phi(E)\neq1$ because $\phi(E)(1-a)=0$ and $a\neq1$. Thus, $\phi(E)=0$, and hence $E\in\ker(\phi|_{\mathcal {P}})=J(r(H)).$ It follows that $E=0$.
If $1-E\in\mathcal {P}$, one can show in a similar way that $1-E=0$, and therefore $E=1$. \qed


Note that Theorem \ref{th21} and Theorem \ref{th1.3} hold as well for the Green rings of Hopf algebras of nilpotent type. In \cite[Theorem 5.4]{WLZ} we showed that the Jacobson radical of the Green ring of a finite dimensional pointed rank one Hopf algebra of nilpotent type is a principal ideal. The method we applied there was to compute the dimension of the Jacobson radical of the complexified Green algebra $R(H)=r(H)\otimes_{\mathbb{Z}}\mathbb{C}$ by determining the irreducible representations of $R(H)$. With the help of the Cartan map, the study of the Jacobson radical of the Green ring becomes more handy.

\section{Stable Green rings and Frobenius-Perron dimensions}
Let $A$ be a finite dimensional non-semisimple Hopf algebra over the complex field $\kk:=\mathbb{C}$. Denote by $A$-mod the category of finite dimensional left $A$-module. The \textit{stable category} $A$-\underline{mod} has the same objects as $A$-mod does, and the space of morphisms from $X$ to $Y$ in $A$-\underline{mod} is the quotient space
$$\text{\underline{Hom}}_A(X,Y):=\text{Hom}_A(X,Y)/\mathcal {P}(X,Y)$$
where $\mathcal {P}(X,Y)$ is the subspace of $\text{Hom}_A(X,Y)$ consisting of morphisms factoring through projective modules.

The stable category $A$-\underline{mod} is a triangulated \cite{Hap} monoidal category with the monoidal structure stemming from that of $A$-mod. The Green ring $r_{st}(A)$ of the stable category $A$-\underline{mod} is called the \textit{stable Green
ring} of $A$. Obviously, the stable Green ring $r_{st}(A)$ admits a $\mathbb{Z}$-basis consisting of all isomorphism classes of finite dimensional indecomposable non-projective $A$-modules.

Let $R$ be a unital associative ring which is free as a $\mathbb{Z}$-module with a finite basis $B=\{b_i\mid i\in I\}$ containing 1. Then $R$ is a \emph{fusion ring} if the following conditions hold (cf. \cite[Definition 1.42.2]{etingof}):
\begin{enumerate}
  \item For any $i,j\in I$, $b_ib_j=\sum_kp_{ij}^kb_k$, where $p_{ij}^k\in\mathbb{Z}_{+}$.
  \item There exists a subset $I_0\subset I$ such that $\sum_{i\in I_0}b_i=1$.
  \item Let $\psi:R\rightarrow\mathbb{Z}$ be the group homomorphism defined by $\psi(b_i)=1$ if $i\in I_0$, and 0 otherwise. There exists an involution $i\mapsto i^*$ of $I$ such that the induced map $b=\sum_{i\in I}p_ib_i\mapsto b^*=\sum_{i\in I}p_ib_{i^*}$ for $p_i\in\mathbb{Z}$ is an anti-involution of ring $R$ and such that $\psi(b_ib_j)=1$ if $i=j^*$, and 0 otherwise.
\end{enumerate}
The fusion ring $R$ is \emph{transitive} if for any $b_i$ and $b_j$, there exist $b_k$ and $b_l$ such that $b_ib_k$ and $b_lb_i$ involve $b_j$ with a nonzero coefficient, for $i,j,k,l\in I$ (cf. \cite[Definition 1.45.1]{etingof}).

For the case when $H$ is a finite dimensional pointed rank one Hopf algebra of non-nilpotent type, the stable Green ring $r_{st}(H)$ of $H$ coincides with the stable Green ring $r_{st}(\overline{H})$ of $\overline{H}$, where the latter has been studied in \cite[Proposition 6.1]{WLZ}. Accordingly, the stable Green ring $r_{st}(H)$ is isomorphic to the quotient ring $r(H)/\mathcal {P}$. More precisely, it is isomorphic to the quotient of the polynomial ring $r(\kk \overline{G})[z]$ modulo the ideal generated by $F_n(a,z)$, where $F_n$ is the Dickson polynomial (of the second type) defined recursively as
$F_1(Y,Z)=1,$ $F_2(Y,Z)=Z$ and $F_j(Y,Z)=ZF_{j-1}(Y,Z)-YF_{j-2}(Y,Z),$ for $j\geq3$.

\begin{prop}
The stable Green ring $r_{st}(H)$ is a transitive fusion ring.
\end{prop}
\proof Let $I=\{(j,i)\mid i\in\Omega_0,1\leq j\leq n-1\}$. Then the set $\{M[j,i]\mid (j,i)\in I\}$ consists of all non-projective indecomposable $H$-modules. It follows that $\{\overline{M[j,i]}\mid (j,i)\in I\}$ forms a basis of $r_{st}(H)=r(H)/\mathcal {P}$. Note that the duality $*$ operator on $r(H)$ is given by $M[j,i]^*=M[j,\tau^{1-j}(i^*)]$ (cf. Proposition \ref{3.7} (1)). This induces an involution on the index set $I$ as $(j,i)^*:=(j,\tau^{1-j}(i^*))$, for any $(j,i)\in I$. Define the group homomorphism $\psi:r_{st}(H)\rightarrow\mathbb{Z}$ by $\psi(\overline{x})=(\delta^*_{[\kk]},x)$, for any $\overline{x}\in r_{st}(H)$. This map is well-defined since $(\delta^*_{[\kk]},x)=0$ if $x\in \mathcal {P}$.
It is straightforward to verify that $r_{st}(H)$ satisfies the definition of a fusion ring given above.
The stable Green ring $r_{st}(H)$ is transitive. We omit the proof because it is similar to the proof of \cite[Proposition 1.45.2]{etingof} using Lemma \ref{19.2} (1).
\qed

The fact  that $r_{st}(H)$ is a transitive fusion ring allow us to calculate the Frobenius-Perron dimensions of objects of the stable category $H$-\underline{mod} in the framework of the stable Green ring $r_{st}(H)$. Let $\text{FPdim}(\overline{M[j,i]})$ be the maximal nonnegative eigenvalue of the matrix of left multiplication by $\overline{M[j,i]}$ with respect to the basis $\{\overline{M[j,i]}\mid (j,i)\in I\}$ of $r_{st}(H)$. Then $\text{FPdim}(\overline{M[j,i]})$ is the \textit{Frobenius-Perron dimension} of the object $M(j,i)$ in $H$-\underline{mod}. The function $\text{FPdim}$ has the following properties, which can be seen directly from \cite[Section 1.45]{etingof}.

\begin{prop}\label{prrr} For any $(j,i)\in I$, the following hold:
\begin{enumerate}
  \item The function $\text{FPdim}:r_{st}(H)\rightarrow\mathbb{C}$ is a ring homomorphism.
  \item $\text{FPdim}$ is a unique nonzero character of $r_{st}(H)$ which takes positive values at each $\overline{M[j,i]}$.
  \item $\text{FPdim}(\overline{M[j,i]^*})=\text{FPdim}(\overline{M[j,i]})$.
 \item $\text{FPdim}(\overline{M[j,i]})\geq1$. If $\text{FPdim}(\overline{M[j,i]})<2$, then $\text{FPdim}(\overline{M[j,i]})=2\cos\frac{\pi}{n}$, for some integer $n>2$.
\end{enumerate}
\end{prop}

Remark that the function $\text{FPdim}$ restricting to the subring $r(\kk \overline{G})$ of $r_{st}(H)$ gives the usual dimensions of $\kk \overline{G}$-modules, namely, $\text{FPdim}([M])=\dim(M)$, for any $\kk \overline{G}$-module $M$, see \cite[Example 1.45.6]{etingof}.

\begin{thm}
$\text{FPdim}(\overline{M[j,i]})=\dim(V_i)F_j(1,2\cos\frac{\pi}{n})$, for any $(j,i)\in I$, where $F_j$ is the $j$-th Dickson polynimial.
\end{thm}
\proof By \cite[Proposition 4.1]{WLZ}, one is able to verify by induction on $j$ that $$F_j(a,M[2,0])=M[j,0],\ \text{for}\ 1\leq j\leq n.$$
Note that $\overline{M[n,0]}=0$ holds in $r_{st}(H)$ since $M(n,0)$ is a projective $H$-module. Consequently, we have:
\begin{align*}0&=\text{FPdim}(\overline{M[n,0]})\\
&=\text{FPdim}(F_n(\overline{a},\overline{M[2,0]}))\\
&=F_n(\text{FPdim}(\overline{a}),\text{FPdim}(\overline{M[2,0]}))\\
&=F_n(1,\text{FPdim}(\overline{M[2,0]})).
\end{align*}
Note that the equality $F_n(1,Z)=0$ has $n-1$ distinct roots: $2\cos\frac{k\pi}{n}$, for $1\leq k\leq n-1$ (cf. \cite[Lemma 5.1]{WLZ}). It follows that $\text{FPdim}(\overline{M[2,0]})$ is equal to $2\cos\frac{k\pi}{n}$ for some $1\leq k\leq n-1$ satisfying $k\mid n$ by Proposition \ref{prrr} (4). We claim that $k=1$. Otherwise, $n_1:=\frac{n}{k}<n$ and $\text{FPdim}(\overline{M[n_1,0]})=F_{n_1}(1,2\cos\frac{\pi}{n_1})=0,$ a contradiction to Proposition \ref{prrr} (4). 
Now for any $(j,i)\in I$, we have
\begin{align*}\text{FPdim}(\overline{M[j,i]})&=\text{FPdim}(\overline{[V_i]})\text{FPdim}(\overline{M[j,0]})\\
&=\dim(V_i)\text{FPdim}(F_j(\overline{a},\overline{M[2,0]}))\\
&=\dim(V_i)F_j(\text{FPdim}(\overline{a}),\text{FPdim}(\overline{M[2,0]}))\\
&=\dim(V_i)F_j(1,2\cos\frac{\pi}{n}).
\end{align*} The proof is completed. \qed

Extending $\text{FPdim}$ from the basis of $r_{st}(H)$ to $R_{st}(H):=\mathbb{C}\otimes_{\mathbb{Z}}r_{st}(H)$ by linearity, we obtain a function $\text{FPdim}$ from $R_{st}(H)$ to $\mathbb{C}$. Denote by $b_{(j,i)}:=\text{FPdim}(\overline{M[j,i]})\overline{M[j,i]}$, for any $(j,i)\in I$. Then the set $\mathbf{B}=\{b_{(j,i)}\mid (j,i)\in I\}$ is a basis of $R_{st}(H)$. The stable Green algebra $R_{st}(H)$ admits some special properties, e.g., the quadruple $(R_{st}(H),\text{FPdim},\mathbf{B},*)$ is a group-like algebra, and hence a bi-Frobenius algebra, see \cite[Section 6]{WLZ} for details.

\section{The Green rings of the Radford Hopf algebras}
In this section, we apply the results
obtained in the previous sections to a family of Hopf algebras, known as
the Radford Hopf algebras. The Radford Hopf algebras was introduced by
Radford in \cite{Rad} so as to give an example of Hopf algebras whose
Jacobson radical is not a Hopf ideal. As we shall see that this family
of Hopf algebras are pointed of rank one and can be derived from
group data of non-nilpotent type.

Let $G$ be a cyclic group of order $mn$ ($m>1$) generated by $g$. Suppose $V_i$ is a one dimensional vector space such that the action of $g$ on $V_i$ is
the scalar multiply by $\omega^i$, where $\omega$ is a primitive
$mn$-th root of unity. Then  $\{V_i\mid i\in\mathbb{Z}_{mn}\}$ forms a complete set of
simple $\mathbbm{k}G$-modules up to isomorphism. Let $\chi$ be the $\mathbbm{k}$-linear
character of $V_{m(n-1)}$. Namely,
$\chi(g)=\omega^{m(n-1)}=\omega^{-m}$, where $\omega^{-m}$ is a
primitive $n$-th root of unity. The order of $\chi$ is $n$ and the $\mathbbm{k}$-linear
character of $V_{m}$ is $\chi^{-1}$.

Let $\mathcal {D}=(G,\chi,g,1)$. Then the group datum $\mathcal {D}$ is of non-nilpotent type since $g^n-1\neq0$ and $\chi^n=1$.
Let $H$ be the Hopf algebra associated to the group
datum $\mathcal {D}$. Then $H$ is generated as an algebra by $g$ and $y$ subject to
$$g^{mn}=1,\ yg=\chi(g)gy=\omega^{-m}gy,\ y^n=g^n-1.$$
The comultiplication $\bigtriangleup$, counit $\varepsilon$, and
antipode $S$ are given respectively by (\ref{3.2}) and (\ref{3.1}).
$H$ is a finite dimensional pointed Hopf algebra of rank one with a
$\mathbbm{k}$-basis $\{y^ig^j\mid 0\leq i\leq n-1,\ 0\leq j\leq mn-1\}$
and $\dim H=mn^2$. The Hopf algebra  $H$ is called the Radford Hopf algebra for the given integers $m$ and $n$.

Let $N=\{1,g^n,g^{2n},\cdots,g^{(m-1)n}\}$,
$\overline{G}=G/N$ and $\overline{\chi}$ the $\mathbbm{k}$-linear
character of $G/N$ such that
$\overline{\chi}(\overline{g^i})=\chi(g^i)$, for $0\leq i\leq mn-1$. Then the Hopf algebra associated to the group datum $\overline{\mathcal
{D}}=(\overline{G},\overline{\chi},\overline{g},0)$ of nilpotent type is nothing but the Taft Hopf algebra $T_n$.
Let $e=\frac{1}{m}\sum_{k=0}^{m-1}g^{kn}$. Then $e$ is a central
idempotent of $H$ and $H$ has the decomposition $H=He\oplus H(1-e)$. By Proposition \ref{propT2}, the quotient $H/H(1-e)$ is a Hopf algebra  isomorphic to $T_n$ and the subalgebra $H(1-e)$ is semisimple.

Denote by $\Omega_0$ the subset of $\mathbb{Z}_{mn}$ consisting of elements divisible by $m$ and $\Omega_1$ the complementary subset of $\Omega_0$.
Let $\tau$ be the permutation of $\mathbb{Z}_{mn}$ determined by
$V_{\chi^{-1}}\otimes V_i\cong V_{\tau(i)}$, where $V_{\chi^{-1}}$
is exactly the simple $\mathbbm{k}G$-module $V_m$ with the character $\chi^{-1}$. It
is easy to see that
$\tau(i)=m+i$, for any $i\in\mathbb{Z}_{mn}$. Let $\langle\tau\rangle$ be the subgroup of the symmetry group $\mathbb{S}_{mn}$
generated by the permutation $\tau$. Then $\langle\tau\rangle$ acts
on the index set $\mathbb{Z}_{mn}$. With this action, the index set $\mathbb{Z}_{mn}$ is divided into $m$
distinct $\langle\tau\rangle$-orbits $[0],[1],[2],\cdots,[m-1],$ where
$[i]=\{i,m+i,2m+i,\cdots,(n-1)m+i\}$, for $0\leq i\leq m-1$. Moreover, $\Omega_0=[0]$ and $\Omega_1=[1]\cup[2]\cup\cdots\cup[m-1]$.

It follows from Theorem \ref{theorem} that $$\{M(k,i),P_{[j]}\mid i\in\Omega_0,1\leq k\leq n,1\leq j\leq m-1\}$$ forms a complete set of finite dimensional indecomposable $H$-modules up to isomorphism. Observe that $V_i\otimes V_j\cong V_{i+j}$ and $(\omega^i\omega^j)^n=1$ if and only if $m\mid i+j$ for any $i,j\in\mathbb{Z}_{mn}$. By
Proposition \ref{propT3} and Proposition \ref{propT4}, we have the following.
\begin{prop}\label{propT6}Let $i,j\in\Omega_1$, $s\in\Omega_0$, $1\leq k\leq n$ and $a=[V_{\chi^{-1}}]=[V_{m}]$.
\begin{enumerate}
\item  $P_{[i]}P_{[j]}=\begin{cases}
(1+a+\cdots+a^{n-1})M[n,0], & m\mid i+j,\\
nP_{[i+j]}, & m\nmid i+j.
\end{cases}$
\item $M[k,s]P_{[j]}=kP_{[j]}$. Moreover $[V_s]P_{[j]}=P_{[j]}$ and $M[2,0]P_{[j]}=2P_{[j]}$.
\end{enumerate}
\end{prop}

Let $\mathbb{Z}[Y,Z,X_1,X_2,\cdots,X_{m-1}]$ be the polynomial ring over $\mathbb{Z}$ in the variables $Y,Z,X_1,X_2,\cdots,X_{m-1}$ and $I$ the ideal of the polynomial ring generated by the following elements:
\begin{equation}\label{equT8}Y^n-1,\ (1+Y-Z)F_n(Y,Z),\ YX_1-X_1,\ ZX_1-2X_1,\end{equation}
\begin{equation}\label{equT15}X_1^j-n^{j-1}X_j,\ \textrm{for}\ 1\leq j\leq m-1,\end{equation}
\begin{equation}\label{equT9}X_1^m-n^{m-2}(1+Y+\cdots+Y^{n-1})F_n(Y,Z).\end{equation}

\begin{thm}\label{thm4.4.2}
The Green ring $r(H)$ of the Radford Hopf algebra $H$ is isomorphic to the quotient ring $\mathbb{Z}[Y,Z,X_1,X_2,\cdots,X_{m-1}]/I$. The Jacobson radical of the quotient ring is a principal ideal generated by $(1-\overline{Y})F_n(\overline{Y},\overline{Z})$.
\end{thm}
\proof Let $r(T_n)$ be the Green ring of the Taft algebra $T_n$. Then $r(T_n)$ is isomorphic to the quotient ring of the polynomial ring $\mathbb{Z}[Y,Z]$ modulo the relations $Y^n=1$ and $(1+Y-Z)F_n(Y,Z)=0$, see \cite[Theorem 3.10]{Chen} or \cite[Theorem 4.3]{WLZ}. Note that $r(H)$ is commutative and is generated as a ring by $P_{[j]}$, for $1\leq j\leq m-1$, over the subring
$r(T_n)$. There is a unique ring
epimorphism $\Phi$ from $\mathbb{Z}[Y,Z,X_1,X_2,\cdots,X_{m-1}]$ to
$r(H)$ such that
$$\Phi(Y)=a,\ \Phi(Z)=M[2,0],\ \Phi(X_j)=P_{[j]},\ \textrm{for}\ 1\leq j\leq m-1.$$
It follows from Proposition \ref{propT6} that the map $\Phi$ vanishes at the generators of the ideal $I$ given by
(\ref{equT8})-(\ref{equT9}). Hence $\Phi$ induces a unique
ring epimorphism $\overline{\Phi}$ from $\mathbb{Z}[Y,Z,X_1,X_2,\cdots,X_{m-1}]/I$ to $r(H)$ such that
$\overline{\Phi}(z+I)=\Phi(z)$, for any $z$ in $\mathbb{Z}[Y,Z,X_1,X_2,\cdots,X_{m-1}]$. Observe
that $\mathbb{Z}[Y,Z,X_1,X_2,\cdots,X_{m-1}]/I$ as a $\mathbb{Z}$-module
has a $\mathbb{Z}$-basis
$\{\overline{Y^iZ^k},\overline{X_j}\mid 0\leq i,k\leq n-1,1\leq j\leq m-1\}$. Thus, as free $\mathbb{Z}$-modules, $\mathbb{Z}[Y,Z,X_1,X_2,\cdots,X_{m-1}]/I$
and $r(H)$ both have the same rank $n^2+m-1$.
As a result, the map $\overline{\Phi}$ is an isomorphism. Thanks to Theorem \ref{th4.3.4}, we have that the Jacobson radical of the quotient ring $\mathbb{Z}[Y,Z,X_1,X_2,\cdots,X_{m-1}]/I$ is a principal ideal generated by $(1-\overline{Y})F_n(\overline{Y},\overline{Z})$.\qed

By Proposition \ref{prop4.50} and Theorem \ref{thm4.4.2}, we have the following corollary.
\begin{cor}
The Grothendieck ring $G_0(H)$ of the Radford Hopf algebra $H$ is isomorphic to the quotient ring $\mathbb{Z}[Y,X_1,X_2,\cdots,X_{m-1}]/I_0$, where $I_0$ is the ideal of  $\mathbb{Z}[Y,X_1,X_2,\cdots,X_{m-1}]$ generated by $Y^n-1,YX_1-X_1,X_1^j-n^{j-1}X_j$ for $1\leq j\leq m-1$ and $X_1^m-n^{m-1}(1+Y+\cdots +Y^{n-1})$.
\end{cor}

\end{document}